\documentclass[10pt]{amsart}

\usepackage{amssymb}
\usepackage[T1]{fontenc} 
\usepackage{wasysym} 
\usepackage{url}
\usepackage{color}

\theoremstyle{plain}
\newtheorem{theorem}{Theorem}[section]
\newtheorem{proposition}[theorem]{Proposition}
\newtheorem{lemma}[theorem]{Lemma}
\newtheorem{corollary}[theorem]{Corollary}
\newtheorem{fact}[theorem]{Fact}

\newtheorem*{claim}{Claim}
\theoremstyle{definition}

\newtheorem{example}[theorem]{Example}
\newtheorem{remark}[theorem]{Remark}

\newtheorem{problem}[theorem]{Problem}

\newcommand{\0}{\emptyset}
\renewcommand{\=}{\approx}
\renewcommand{\leq}{\leqslant}
\renewcommand{\geq}{\geqslant}
\newcommand{\meet}{\wedge}
\newcommand{\bigmeet}{\bigwedge}
\newcommand{\join}{\vee}
\newcommand{\bigjoin}{\bigvee}

\newcommand{\mmeet}{\protect{\wedge\mbox{}\kern-5.2pt\wedge}}
\newcommand{\mbigmeet}{\bigwedge\mbox{}\kern-10pt\bigwedge}
\newcommand{\dsmbigmeet}{\bigwedge\mbox{}\kern-13.3pt\bigwedge} 
\newcommand{\mjoin}{\protect{\vee\mbox{}\kern-5.2pt\vee}}
\newcommand{\mbigjoin}{\bigvee\mbox{}\kern-10pt\bigvee}
\newcommand{\mneg}{\text{\raisebox{-.5pt}{$\neg$}\mbox{}\kern-6.7pt\raisebox{1pt}{$\neg$}}}
\newcommand{\mto}{\Rightarrow}
\newcommand{\mBox}{\text{\raisebox{0.9pt}{$\Box$}}}

\usepackage[all,2cell,ps]{xy}
\newcommand{\DiamondBox}{\protect{\kern 0.88pt\text{\fontsize{7pt}{12pt}\selectfont\raisebox{1.02pt}{$\Diamond$}}\kern-6.605pt\mBox}}
\newcommand{\piccirc}{\protect{\text{\fontsize{14pt}{12pt}\selectfont$\circ$}}}

\newcommand{\Con}{\operatorname{Con}}
\newcommand{\id}{\operatorname{id}}
\newcommand{\Th}{\operatorname{Th}}

\newcommand{\V}{\mathcal V}
\newcommand{\W}{\mathcal W}
\newcommand{\U}{\mathcal U}
\newcommand{\K}{\mathcal K}
\newcommand{\Q}{\mathcal Q}
\newcommand{\dS}{\mathcal S}
\renewcommand{\S}{\mathbf S}
\renewcommand{\P}{\mathbf P}
\newcommand{\M}{\mathbf{M}}
\newcommand{\N}{\mathbf{N}}
\newcommand{\A}{\mathbf{A}}
\newcommand{\B}{\mathbf{B}}
\newcommand{\C}{\mathbf C}
\newcommand{\R}{\mathbf R}
\renewcommand{\1}{\mathbf{1}}
\renewcommand{\2}{\mathbf{2}}
\newcommand{\4}{\mathbf{4}}

\newcommand{\Fr}{\mathbf{F}}
\newcommand{\G}{\mathbf{G}}
\renewcommand{\H}{\mathbf{H}}
\newcommand{\Lev}{\mathbf{Lev}}
\newcommand{\Nat}{\mathbb N}

\newcommand{\Form}{\mathbf{Form}}

\begin{document}

\title[Almost structural completeness]{Almost structural completeness;\\ an algebraic approach}

\author{Wojciech Dzik}

\address{Institute of Mathematics,
University of Silesia,
ul. Bankowa 14,
40-007 Katowice, Poland}
\email{dzikw@ux2.math.us.edu.pl}

\author{Micha\l{}~M. Stronkowski}

\address{Faculty of Mathematics and Information Sciences,
Warsaw University of Technology, ul. Koszykowa 75, 00-662
Warsaw, Poland}
\email{m.stronkowski@mini.pw.edu.pl}

\thanks{The wok of the second author was supported by the Polish National Science Centre grant no. DEC- 2011/01/D/ST1/06136.}

\keywords{Almost structural completeness, structural completeness, quasivarieties, axiomatization, modal normal logics, varieties of closure algebras, equationally definable principal relative congruences, finite model property.}

\subjclass[2010]{08C15, 03G27, 03B45 , 03B22, 06E25.}


\begin{abstract} A deductive system is structurally complete if its admissible inference rules
are derivable. 
For several important systems, like modal logic S5, failure of structural completeness is caused only by the underivability of passive rules, i.e. rules that can not be applied to theorems of the system. 
Neglecting passive rules leads to the notion of almost structural completeness, that means, derivablity of admissible non-passive rules. 
 Almost structural completeness for quasivarieties and varieties of general algebras is investigated here by purely algebraic means. The results apply to all algebraizable deductive systems.

Firstly, various characterizations of almost structurally complete quasivarieties are presented. Two of them are general:
expressed with finitely presented algebras, and with subdirectly irreducible algebras. One is restricted to quasivarieties with finite model property and equationally definable principal relative congruences, where the condition is verifiable on finite subdirectly irreducible algebras.

Secondly, examples of almost structurally complete varieties are provided
Particular emphasis is put on varieties of closure algebras, that are known to constitute adequate semantics for normal extensions of S4 modal logic.
A certain infinite family of such almost structurally complete, but not structurally complete, varieties is constructed.
Every variety from this family has a finitely presented unifiable algebra which does not embed into any free algebra for this variety. Hence unification in it is not unitary. This shows that almost structural completeness is strictly weaker than  projective unification for varieties of closure algebras.
\end{abstract}

\maketitle

\section{Introduction}

In order to present motivation for the paper and for almost structural completeness let us recall basic notions from algebraic logic.
Let $\mathcal L$ be a propositional language, i.e., a set of logical connectives with ascribed arities, and let $\Form$ be the algebra of formulas in $\mathcal L$ over a denumerable set of variables.
An {\emph{inference}) \emph{rule} is a pair from $\mathcal P(Form)\times Form$, written as $\Phi\slash\varphi$, where $\mathcal P(Form)$ is the powerset of $Form$. By a \emph{deductive system} we mean a pair $\mathcal S=(\Form,\vdash)$, where $\vdash$ is a finitary structural (i.e. preserving substitutions) consequence relation, this is a set of rules satisfying  appropriate postulates\footnote{We adopt the definition form
\cite{FJP03}. However it is also a common practice to use the term ``deductive system'' for the basis of deductive system in our sense.} \cite{Cze01,FR09,FJP03,PW08,Raf11b,Ras74,Woj73,Woj98}. (We drop here most of definitions and keep such a level of formality that, allows the reader to comprehend the main ideas.)
Let $\Th(\dS)=\{\varphi\in Form\mid \0\vdash \varphi\}$ be the set of \emph{theorems} of $\dS$. A \emph{basis} or an \emph{axiomatization} of $\dS$ is a pair $(A,R)$, where $A\subseteq\Th(\dS)$ and $R\subseteq\;\vdash$ are such that $\Phi\vdash\varphi$ iff there is a \emph{proof (derivation)} form $A\cup \Phi$ for $\varphi$ by means of rules from $R$.

Often, instead of a deductive system $\dS$,  interest is put mainly on  the set of its theorems $L=\Th(\dS)$.  The set $L$ is then called  a \emph{logic}. It happens especially when $R$ is chosen in some default way. For instance, for intermediate logics $R$ consists of \emph{Modus Ponens} and for normal modal logics $R$ consists of \emph{Modus Ponens} and \emph{Necessitation rule}. Given a basis $A$ of a logic $L$, equipped with a default set of rules $R$, a formula can be proved or derived from $A$. Proofs of theorems may be shorten by allowing new rules. Such extension of $R$ may be done in two ways:
\begin{enumerate}
\item by adding \emph{derivable} rules, i.e., those that are in $\vdash$,
\item by adding \emph{admissible} but non-derivable rules, i.e., under which the set of theorems is closed but which are not derivable.
\end{enumerate}
The admissibility is more elusive than the derivability. Its verification for a rule may be a challenging task \cite{Ryb97}. Deductive systems (and logics) for which all admissible rules are derivable are called \emph{structurally complete} (SC for short).


A rules $\Phi\slash\varphi$ is \emph{passive} if for every substitution $\sigma$  (i.e., an endomorphism of $\Form$)  the set $\sigma(\Phi)$ is not contained in $\Th(\dS)$. Such a rule can not be applied to theorems. 
There are important examples of deductive systems that are not SC and in which admissible non-derivable rules are passive. Such systems are called \emph{almost structurally complete} (ASC for short).
In this case a proof of any theorem  also cannot shortened  by the method (2). 

In the following example of modal logic an advantage of ASC over SC is particularly apparent. 
Let $L$ be a modal normal logic with a basis $A$. Recall that $L$ has an adequate algebraic semantics given by a variety $\V$ of modal algebras (see Section \ref{sec:: examples} for definitions of modal and closure algebras). A formula $\varphi(\bar x)$ holds in a modal algebra $\M$ provided \mbox{$\M\models(\forall\bar x)\,\varphi(\bar x)\= 1$}, and, a rule $\varphi_1(\bar x),\ldots,\varphi_n(\bar x)\slash\varphi(\bar x)$ (we adopt a common convention and drop the curly brackets) holds in $\M$ if the quasi-identity $(\forall \bar s)[\varphi_1(\bar x)\=1\meet\cdots\meet\varphi_n(\bar x)\=1\to\varphi(\bar x)\=1]$ holds in $\M$. Then a formula belongs to $L$ iff it holds in all algebras from $\V$, and similarly a rule is derivable iff it holds in all algebras from $\V$.

Assume that algebras $\2$ and $\S_2$, depicted in Figure \ref{fig:: 2 and S_2}, belong to $\V$.
\begin{figure}[h]
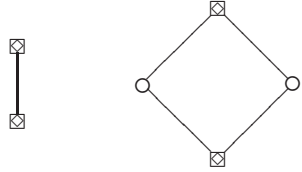

\[
\xy
(10,0)*{}; 
(10,15)*{\DiamondBox}="AC";
(10,5)*{\DiamondBox}="ZERO";
(10,6.1); (10,14.2) **\dir{-}; 
\endxy
\quad\quad\quad\quad
\xy
(10,20)*{\DiamondBox}="AC";
(0,9.8)*{\piccirc}="A";  (20,10)*{\piccirc}="C";
(10,0)*{\DiamondBox}="ZERO";
(9.07,1.07); (.7,9.44) **\dir{-}; 
(10.93,1.07);  (19.35,9.49) **\dir{-}; 
(9.07,19.22); (.55,10.7) **\dir{-}; 
(19.32,10.82); (10.9,19.24) **\dir{-};  
\endxy
\]
\caption{The modal algebras $\bf 2$ and ${\bf S}_2$.}
\label{fig:: 2 and S_2}
\end{figure}
The algebras $\2$ and $\S_2$ are closure algebras with only the top and the bottom elements closed (and open).
Then $L$, which corresponds to $\V$,   cannot be SC. Indeed, the rule $r$
\[
\Diamond p,\Diamond\neg p\;\,\slash\perp.
\]
is not valid in $\S_2$. Therefore it is not derivable. However, $r$ is admissible, because
$\2\models(\forall y)[\Diamond y \not\=1 \,\join\,\Diamond\neg y\not\=1]$. Hence \mbox{$\2\not\models(\forall \bar x)[\Diamond \varphi(\bar x)\=1 \,\meet\,\Diamond\neg \varphi(\bar x)\=1]$} and, consequently, $\{\Diamond\varphi(\bar x) ,\Diamond\neg \varphi(\bar x)\}\not\subseteq L$ for every modal formula $\varphi(\bar x)$. Thus $r$ is passive and, therefore, admissible.

There are many normal modal logics of this kind which are ASC. One of them is modal logic S5 (the logic of equivalence relations). An algebraic semantics for S5 is given by the variety of monadic algebras,  at which we will  look in  Example \ref{exp:: monadic algebras}.


The notion of ASC  
will be  studied here by means of algebraic semantics.
Recall that every algebraizable deductive system $\mathcal S$ has an adequate semantics which is a quasivariety $\Q$ of algebras \cite{BP89}. In particular: logical connectives become basic operations, formulas become terms, theorems of $\mathcal S$ correspond to identities true in $\Q$, and derivable rules of $\mathcal S$ correspond to quasi-identities true in $\Q$. Therefore the property of being ASC may be formulated for quasivarieties.

In the first part of the paper (Sections \ref{sec:: general characterization}-\ref{sec: criterion}) we develop a general theory of ASC for quasivarieties. We present various characterizations of ASC quasivarieties.
Two of them are general:  expressed with finitely presented algebras (Theorem \ref{thm::main characterization} Condition (5)), and with subdirectly irreducible algebras (Theorem \ref{thm::main characterization} Condition (3)). One is restricted to quasivarieties with finite model property and equationally definable principal relative congruences where the condition is verifiable on finite subdirectly irreducible algebras (Theorem \ref{thm::criterion}). We would like to note that the condition for being ASC with finitely presented algebras has connection to unification theory, see Section \ref{sec:: projective unification}). But the condition with subdirectly irreducible algebras is purely algebraic and probably could not be discovered without algebraic tools.

In the second part of the paper (Sections \ref{sec:: examples} and \ref{sec:: more examples}) we illustrate theoretical considerations by showing how our results  may be used to establish ASC for particular varieties. Until now the most common method for proving ASC was by an application of projective unification see e.g. \cite{Dzi08,Dzi11,DW12,DW14}, (see \cite{MR13} for an exception). In turn, we use algebraic tools like subdirectly irreducible algebras and free algebras.

We put an emphasis on varieties of closure algebras. They  constitute adequate semantics for transitive reflexive normal modal logics, i.e., for normal extensions of S4 modal logic. The main result here is the description of an infinite family of ASC, but not SC, varieties without unitary unification and with finitely presented unifiable algebras not embeddable into free algebras (Theorem \ref{thm:: modal companions of Medvedev and Levin}). Thus a verification of ASC for them could not be obtained by means of projective unification. The idea behind the construction is the following one. We consider a variety $\U$ of SC closure algebras without projective unification. We ``spoil it a bit'' by taking the varietal join $\U\join\W$ with a non-minimal variety $\W$ of monadic algebras. Non-minimal varieties of monadic algebras are known to be ASC but not SC. In order to prove that such join is still ASC we have to develop the theory of closure algebras. In particular, on the way, we show that an ASC variety of closure algebras is SC iff it satisfies McKinsey identity (Proposition \ref{prop:: SC <-> McK}). Moreover, we describe free algebras for $\U\join\W$ by means of free algebras for $\U$ and $\W$, where $\U$ is a variety of McKinsey algebras and $\W$ is a variety of monadic algebras (Proposition \ref{prop:: free algebras in McKinsey join monadic}). Finally, for $\U$ we can take the smallest modal companions of Levin and Medvedev varieties of Heyting algebras  described in Example \ref{exp:: Levin and Medvedev}. These varieties characterize Levin or Medvedev intermediate logics which are known to be SC \cite{Pruc76} and not possessing projective unification \cite{Dzi06b}.

The results may be applied to the axiomatization problem for quasivarieties and to find a basis of admissible rules of a deductive systems. Indeed, if ASC property is established for a quasivariety $\Q$ (a deductive system $\dS$), then an axiomatization of every of its subquasivariety containing all free algebras for $\Q$ (the extension of $\dS$ obtained by adding some admissible rules) may be obtained by adding passive quasi-identities (rules). In particular, an axiomatization of the quasivariety generated by all free algebras for $\Q$ (the extension of $\dS$ obtained by adding all admissible rules) may be obtained by adding all passive quasi-identities (all passive rules),
see Section \ref{sec: asc core}. This fact was used in \cite{DW14}, where the analysis of passive rules lead to a description of the lattice of all deductive systems extending modal logic $S4.3$.

\subsection*{Historical notes}

The notion of SC was introduced by Pogorzelski in \cite{Pog71}
and then investigated by many authors. The reader may consult the monograph \cite{PW08} and references therein for older results concerning SC property.
Let us recall here more recent works:  \cite{OR07} about varieties of positive Sugihara monoids, \cite{ORV08}  about substructural logics, \cite{CM09} about fuzzy logics, \cite{Slo12} about some fragment of the intuitionistic logic, \cite{Kow14} about BCK logic, \cite{CSV14} for semisimple varieties and discriminator varieties, and \cite{Raf14} which contains a  general considerations from abstract algebraic logic perspective and results for some non-algebraizable deductive systems.
Although SC property was often investigated algebraically, there are only few papers about it for algebras not connected to logic. Among them the paper of Bergman \cite{Ber88} deserve special attention. In particular, he formulated the condition for a quasivariety $\Q$ to be SC. It says that $\Q$ must be generated by its free algebras. For specific algebras, SC was investigated for lattices in \cite{Igo74} and for modules in \cite{HB85}.






ASC appeared for the first time, though under different name in \cite{Dzi06} and  \cite{Dzi08} as an application of projective unification.
Projective unification and thus ASC was established for some varieties and logics. Probably the most prominent examples are discriminator varieties. This includes varieties of e.g. Boolean algebras, monadic algebras, rings satisfying $x^m\=x$  for a finite $m>1$, MV${}_n$-algebras, $n$-valued Post algebras, cylindric algebras of dimension $n$, all for finite $n$ \cite[Theorem 3.1]{Bur92}. For intuitionistic logic it was shown that every extension of G\"odel-Dummmett logic LC has projective unification \cite{Wro95}. For a normal extension $L$ of S4 modal logics it was proved that $L$ has projective unification iff $L$ is an extension of S4.3 \cite[Corollary 3.19]{DW12}, see also \cite[Theorem 5]{Iem14}. Projective unification was also verified for some others modal logics not extending S4
\cite{Dzi06}, for $k$-potent extensions of basic fuzzy logic and hoops \cite{Dzi08}, and for some Fregean varieties \cite{Slo12b}.

A general investigations on ASC property was, independently from ours, undertaken in \cite{MR13}. In particular, Corollary \ref{cor:: about C} was published there for the first time. Note however that, contrary to our paper, \cite{MR13} is focused on a finitely generated case. The main result there concerning ASC property is the algorithm (with an applicable software) for deciding whether a given finite family of finite algebras (of a small size) in a finite language generates an ASC quasivariety.



Several variants of ``completeness'' for deductive systems other than SC and ASC were proposed, like maximality and Post completeness. 
Passive (or overflow) structural completeness \cite{Wro09} may be considered as complementary to ASC. 
A deductive system is \emph{passively structurally complete} (PSC) if every its admissible and passive rule is derivable. Clearly, a deductive system is SC iff it is ASC and PSC. PSC property was also investigated in the context of fuzzy logics \cite{CM09}.


Finally let us note that our reaserch belongs into an intensively investigated area of  admissibility of rules in general, see for instance, \cite{Iem14b,Jer05,Ryb97}.

\section{Concepts from quasivariety theory}\label{sec::quasivariety theory}
Though most deductive systems we are interested in have algebraic semantics given by varieties (they are strongly algebraizable), the right language to deal with SC, ASC and admissability is quasivariety theory. This is so because we have to work with quasi-identities anyway. Therefore we will formulate our main results for quasivarieties. Let us here recall needed notions and facts from this theory.

Following \cite{BS81, Gor98,Mal70} we call a first order sentence a \emph{quasi-identity} if it is of the form
\[
(\forall \bar x)\, [s_1(\bar x)\=t_1(\bar x)\meet\cdots\meet s_{n}(\bar x)\=t_{n}(\bar x)\to s(\bar x)\=t(\bar x)],
\]
where $n\in\mathbb N$. We allow $n$ to be zero, and in such case we call the sentence an \emph{identity}. It will be convenient to have a more compact notation for quasi-identities, and we will often write them in the form
\[
(\forall \bar x)\, [\varphi(\bar x)\to \psi(\bar x)],
\]
where $\varphi$ is a conjunction of equations (i.e., atomic formulas) and $\psi$ is an equation. We call $\varphi$ the \emph{premise} and $\psi$ the \emph{conclusion} of a quasi-identity.

By a \emph{(quasi-)equational theory} of a class $\K$ of algebras in the same language we mean the set of (quasi-)identities true in $\K$.
A \emph{(quasi)variety} is a class defined by (quasi-)identities. Equivalently, a class of algebras in the same language is a quasivariety if it is closed under taking substructures, direct products and ultraproducts. If it is additionally closed under taking homomorphic images, it is a variety. (We tacitly assume that all considered classes contain algebras in the same language and are closed under taking isomorphic images. Also all considered class operators are assumed to be composed with isomorphic image class operator.) A (quasi)variety is \emph{trivial} if it consists of one-element algebras, and is \emph{minimal} if it properly contains only a trivial (quasi)variety. We say that a class is a (quasi)variety generated by a class $\K$ if it is the smallest (quasi)variety containing $\K$, i.e., the class defined by the \mbox{(quasi-)equational} theory of $\K$. We denote such class by ${\sf V}(\K)$ (${\sf Q}(\K)$ respectively). In case when $\K=\{\A\}$ we simplify the notation by writing ${\sf V}(\A)$ (${\sf Q}(\A)$ respectively). Note that ${\sf V}(\K)={\sf HSP}(\K)$ and ${\sf Q}(\K)={\sf SPP_U}(\K)$, where $\sf H$, $\sf S$, $\sf P$, $\sf P_U$
are homomorphic image, subalgebra, direct product and ultraprodut class operators \cite[Theorems II.9.5 and V.2.25]{BS81}.

Let $\Q$ be a quasivariety. A congruence
$\alpha$ on an algebra $\A$ is called a \mbox{\emph{$\Q$-congruence}}
provided $\A/\alpha\in\Q$. Note that $\A\in\Q$ if and only if the equality relation on $A$ is a $\Q$-congruence. The set $\Con_\Q(\A)$ of all
\mbox{$\Q$-congruences} of $\A$ forms an algebraic lattice which
is a meet-subsemilattice of $\Con(\A)$ of all congruences of $\A$ \mbox{\cite[Corollary 1.4.11]{Gor98}}. We say that a (quasi)variety $\Q$ is (relatively) congruence-distributive if all lattices $\Con_\Q(\A)$ are distributive.

A nontrivial algebra $\S$ is $\Q$-simple if $\Con_\Q(\S)$ has exactly two elements: the equality relation $\id_S$ on $S$ and the total relation $S^2$ on $S$. A nontrivial algebra $\S\in\Q$ is
$\Q$-\emph{subdirectly irreducible} if the equality relation on $A$ is completely meet irreducible in $\Con_\Q(\A)$. (In case when $\Q$ is a variety we do drop the prefix \mbox{``$\Q$-''.)} Let us denote the class of all $\Q$-subdirectly irreducible algebras by $\Q_{SI}$. The importance of $\Q_{SI}$ follows from the fact that this class determines $\Q$. Indeed, in an algebraic lattice each element
is a meet of completely meet-irreducible elements. Moreover, for $\A\in\Q$
the lattice $\Con_\Q(\A)$ is algebraic. Thus we have the following fact.

\begin{proposition}[\protect{\cite[Theorem 3.1.1]{Gor98}}]\label{prop::Mal'cev thm}
Every algebra in a quasivariety $\Q$ is isomorphic to a subdirect product of $\Q$-subdirectly irreducible algebras. In particular, $\Q$ is generated by $\Q_{SI}$.
\end{proposition}

Let $\G\in\Q$ and $X\subseteq G$. We say that $\G$ is \emph{free for $\Q$ over $X$}, and is \emph{of rank} $|X|$, if $\G\in\Q$ and it satisfies the following \emph{universal mapping property}: every mapping $f\colon X\to A$, where $A$ is a carrier of an algebra $\A$ in $\Q$, is uniquely extendable to a homomorphism $\bar f\colon\G\to\A$. Elements of $X$ are then called \emph{free generators of} $\G$.
If $\Q$ contains a nontrivial algebra, then it has free algebras over arbitrary non-empty sets and, in fact, they coincide with free algebras for the variety $\sf V(\Q)$. (Note here that $\sf V(\Q)$ is the class of all homomorphic images of algebras from $\Q$.) Let us fix a countably infinite set of variables $V=\{v_0,v_1,v_2,\ldots\}$. We denote a free algebra for $\Q$ over $V$ by $\Fr$ and the free algebra for $\Q$ over $V_k=\{v_0,v_1,\ldots,v_{k-1}\}$ by $\Fr(k)$. One may construct $\Fr$ and $\Fr(k)$ by taking the algebra of terms over $V$, or $V_k$ respectively, and divide it by the congruence identifying terms $s(\bar v),t(\bar v)$ which determine the same term operation on every algebra from $\Q$ (in other words, when $\Q\models(\forall\bar x)[t(\bar x)\=s(\bar x)]$). The algebra $\Fr$ is an union of a chain of subalgebras which are isomorphic to $\Fr(k)$. It follows that the family of all free algebras for $\Q$ of finite rank generates the quasivariety ${\sf Q}(\Fr)$. We will notationally identify terms with elements of $\Fr$ that they represent.

For an algebra $\A$ and a set $H\subseteq A^2$ there exists the least $\Q$-con\-gruence $\theta_\Q(H)$ on $\A$ containing $H$. When $H=\{(a,b)\}$ we just write $\theta_\Q(a,b)$. (When $\Q$ is a variety we also simplify the notation by dropping the subscript $\Q$.)
We say that an algebra is \emph{$\Q$-finitely presented}
if it is isomorphic the $\Fr(k)/\theta_\Q(H)$ for some natural number $k$ and some finite set $H$ \cite[Chapter 2]{Gor98}. The class of all $\Q$-finitely presented algebras will be denoted by $\Q_{FP}$.
For a tuple $\bar x=(x_0,\ldots,x_{k-1})$ of variables and a conjunction of equations $\varphi (\bar x)=s_1(\bar x)\=t_1(\bar x)\meet\cdots\meet s_n(\bar x)\=t_n(\bar x)$ let
\[
\P_{\varphi(\bar x)}=\Fr(k)/\theta_\Q(\{(s_1(\bar v),t_1(\bar v)),\ldots,(s_n(\bar v),t_n(\bar v))\}),
\]
where $\bar v=(v_0,\ldots,v_{k-1})$. Note that every finitely presented algebra is isomorphic to some $\P_{\varphi(\bar x)}$.
Observe  that $\Q$ satisfies a quasi-identity $(\forall \bar x)\, [\varphi(\bar x)\to \psi(\bar x)]$ iff $\P_{\varphi(\bar x)}\models \psi(\bar v)$ (we notationally identify variables from $\bar v$ with their congruence classes). In particular,  $\Q$ satisfies an identity $(\forall \bar x)\, [s(\bar x)\=t(\bar x)]$ iff $s(\bar v)=t(\bar v)$ in $\Fr(k)$. We say that a $\Q$-finitely presented algebra $\P$  is \emph{unifiable} (\emph{in $\Q$}) provided that there exists a homomorphism from $\P$ into $\Fr$. Every such homomorphism is called an \emph{unifier for $\P$}. Finally, let as recall that $\Q_{FP}$ also generates $\Q$. Strictly we have the following fact.

\begin{proposition}[\protect{\cite[Proposition 2.1.18]{Gor98}}]\label{prop::direct limit thm}
Every algebra in a quasivariety $\Q$ is isomorphic to a direct limit of $\Q$-finitely presented algebras.
\end{proposition}

We will use the following folklore fact.

\begin{fact}\label{fact::finite->fp}
Let $\Q$ be a quasivariety in a finite language. Then every finite algebra in $\Q$ is $\Q$-finitely presented. 
\end{fact} 

\begin{proof}
Let $\P$ be a finite algebra in $\Q$. Take a tuple $\bar a=(a_0,\ldots,a_{k-1})$ such that $P=\{a_0,\ldots,a_{k-1}\}$. Let $\bar x=(x_0,\ldots,x_{k-1})$ be the tuple of variables of length $k$. Define the set $\Phi(\bar x)$ consisting of all equations of the form 
$\omega(x_{i_0},\ldots,x_{i_{n-1}})\=x_{i_n}$, where $\omega$ is an $n$-ary basic operation from the language of $\Q$, variables are from $\bar x$, and $\omega(a_{i_0},\ldots,a_{i_{n-1}})\=a_{i_n}$ in $\P$. Since $P$ is finite and the number of basic operations is finite, $\Phi(\bar x)$ is finite. Let $\varphi(\bar x)=\bigmeet\Phi(\bar x)$. Then $\P\cong\P_{\varphi(\bar x)}$. In order to see this let us consider a homomorphism $h\colon\Fr(k)\to\P$ satisfying $h(\bar x)=\bar a$. Its existence follows from the universal mapping property. Now, since the kernel of $h$ is a $\Q$-congruence and $\P\models \varphi(\bar a)$, the definition of $\P_{\varphi(\bar x)}$ yields that $h$ factors through $g\colon \P_{\varphi(\bar x)}\to\P$, where $g(\bar v)=\bar a$. Clearly, $g$ is surjective. Moreover, the definition of $\varphi(\bar x)$ implies that $P_{\varphi(\bar x)}$ has at most $k$ elements. Thus $g$ must be also injective.

Actually, the proof may be slightly modified in order to obtain a stronger fact: If $l$ is the cardinality of a smallest generating set for $\P$, then in a defining formula $\varphi(\bar x)$ we need only $l$ variables.
\end{proof}

\section{General Characterizations and first observations}\label{sec:: general characterization}

For a quasi-identity $q=(\forall \bar x)\, [\varphi(\bar x)\to \psi(\bar x)]$ let
\[
q^*=(\forall \bar x)\, [\neg\varphi(\bar x)].
\]
We partition the set of quasi-identities true in $\Fr$ into two sets: the set of \emph{$\Q$-active} quasi-identities $q$  for which $q^*$ does not hold in $\Fr$, and the set of \emph{$\Q$-passive} quasi-identities $q$  for which $q^*$  holds in $\Fr$. Equivalently, a quasi-identity $q$ true in $\Fr$ is $\Q$-active if $\P_{\varphi(\bar x)}$ is unifiable and it is $\Q$-passive if $\P_{\varphi(\bar x)}$ is not unifiable, where $\varphi(\bar x)$ is the premise of $q$.


A quasivariety $\Q$ is \emph{structurally complete} (SC for short) provided that every quasi-identity which is true in $\Fr$ is also true in $\Q$, in other words if $\Q={\sf Q}(\Fr)$.
A quasivariety $\Q$ is \emph{almost structurally complete} (ASC for short\footnote{Maybe a better full form of ASC would be \emph{active structural completeness} as Alexander Cytkin privately suggested to us.}) provided that every $\Q$-active quasi-identity holds in $\Q$. We will also use the abbreviation ASC$\setminus$SC to indicate that a considered quasivariety is ASC but is not SC.

Let us start considerations by providing various conditions for quasivarieties equivalent to being ASC.

We will write $\A\to\B$ to code the supposition that there is a homomorphism from $\A$ into $\B$. In particular, for a $\Q$-finitely presented algebra $\P$, $\P\to\Fr$ means that $\P$ is unifiable.

\begin{theorem}\label{thm::main characterization}
The following conditions are equivalent:
\begin{enumerate}
\item $\Q$ is ASC;
\item For every $\A\in\Q$, $\A\times\Fr\in\sf Q(\Fr)$;
\item For every $\S\in\Q_{SI}$, $\S\times\Fr\in\sf Q(\Fr)$;
\item For every $\A\in\Q$, $\A\to\Fr$ yields $\A\in\sf Q(\Fr)$;
\item For every $\P\in\Q_{FP}$, $\P\to\Fr$ yields $\P\in\sf Q(\Fr)$.
\end{enumerate}
\end{theorem}

\begin{proof}
The implications (2)$\Rightarrow$(3) and (4)$\Rightarrow$(5) are obvious. \\
(1)$\Rightarrow$(2) Let $\A\in\Q$ and consider a quasi-identity $q$ true in $\Fr$. We wish to show that $\A\times\Fr\models q$. If $\Q\models q$, then it clearly holds since $\A\times\Fr\in\Q$. So suppose that $\Q\not\models q$. Then, by the definition of ASC, $\Fr\models q^*$. Thus $\A\times\Fr\models q^*$, and therefore $\A\times\Fr\models q$.

\noindent(2)$\Rightarrow$(1) Let $q=(\forall \bar x)\, [\varphi(\bar x)\to \psi(\bar x)]$ and assume that $\Fr\models q$ and $\Q\not\models q$. Then $q$ is not valid in some $\A\in\Q$, i.e., there is a tuple $\bar a$ of elements in $A$ such that $\A\models \varphi(\bar a)\meet\neg\psi(\bar a)$. We would like to show that $q$ is $\Q$-passive. Suppose that, on the contrary, $\Fr\not\models q^*$. This means that there is a tuple $\bar t$ from $F$ such that $\Fr\models\varphi(\bar t)$. Then $\A\times\Fr\models\varphi(\bar d)$, where $\bar d$ is the tuple of pairs of elements from $\bar a$ and $\bar t$ in the respective order. By (2), $\A\times\Fr\models q$, and hence $\A\times\Fr\models\psi(\bar d)$. This yields that $\A\models \psi(\bar a)$, and we obtained a contradiction.

\noindent(3)$\Rightarrow$(2) Let $\A\in\Q$. By Proposition \ref{prop::Mal'cev thm}, $\A$ is isomorphic to a subdirect product of $\S_i\in\Q_{SI}$, $i\in I$. If $I=\0$, then $\A$ is trivial and $\A\times \Fr\cong\Fr\in{\sf Q(\Fr)}$. So let us assume that $I\neq\0$. Then $\Fr$ is isomorphic with the diagonal of $\Fr^I$, and hence $\A\times\Fr$ is isomorphic with a subalgebra of $\A\times \Fr^I$. Further, the latter is isomorphic to a subalgebra of $\left(\prod_{i\in I}\S_i\right)\times\Fr^I\cong \prod_{i\in I}(\S_i\times\Fr)$. Thus (3) yields  that $\A\times\Fr\in{\sf Q}(\Fr)$.

\noindent(2)$\Rightarrow$(4) Assume that there is a homomorphism $h$ from $\A\in\Q$ into $\Fr$. Let $\mathbf{h}$ be the subalgebra of $\A\times\Fr$ with the carrier $h$. By (2) the algebra  $\mathbf h$ belongs to ${\sf Q}(\Fr)$. Since  $\A\cong\mathbf{h}$, the algebra $\A$ also belongs to ${\sf Q}(\Fr)$.

\noindent(4)$\Rightarrow$(2) It holds since there is a homomorphism from $\A\times\Fr$ into $\Fr$, namely, the second projection.

\noindent(5)$\Rightarrow$(4) Let $\A\in\Q$. By Proposition \ref{prop::direct limit thm}, we may assume that $\A$ is a direct limit $\displaystyle{\lim_{\longrightarrow}}\P_i$ of $\Q$-finitely presented algebras $\P_i$. Let $k_i\colon\P_i\to\A$ be the associated canonical homomorphisms. Assume that $f\colon\A\to\Fr$. Then $f\circ k_i\colon\P_i\to\Fr$, and (5) gives $\P_i\in\sf Q(\Fr)$. Since every quasivariety is closed under taking direct limits \cite[Theorem 1.2.12]{Gor98}, $\A$ belongs to $\sf Q(\Fr)$.
\end{proof}

The list of Conditions from Theorem \ref{thm::main characterization} is not full but we consider them as the most fundamental. In this section we will also formulate additional conditions equivalent to ASC which will be used in our considerations.

\begin{corollary}\label{cor:: about C}
Let $\C$ be a subalgebra of $\Fr$, e.g. $\Fr(1)$ or $\Fr(0)$ (if it exists). Then $\Q$ is ASC if and only if one  of the following conditions holds.
\begin{enumerate}
\item[(2')] For every $\A\in\Q$, $\A\times\C\in\sf Q(\Fr)$;
\item[(3')] For every $\S\in\Q_{SI}$, $\S\times\C\in\sf Q(\Fr)$.
\end{enumerate}
\end{corollary}

\begin{proof}
For every algebra $\A$ we have $\A\times\C\leq\A\times\Fr$, and hence conditions (2) and (3) from Theorem \ref{thm::main characterization} yield (2') and (3'), respectively. For proving the converse, let us consider a homomorphism $h\colon\Fr\to\C\leq\Fr$. Its existence is guaranteed by the universal mapping property. Then  $\A\times\Fr$ embeds into $\A\times\C\times\Fr$ via the mapping $(a,t)\mapsto(a,h(t),t)$. This shows that (2') and (3') yields (2) and (3) from Theorem \ref{thm::main characterization}, respectively.
\end{proof}

\begin{remark}
The equivalence of (2') in Corollary \ref{cor:: about C} with ASC was independently proved in \cite[Theorem 18]{MR13}.
\end{remark}

\begin{corollary}\label{cor:: on finitely presented}
A quasivariety $\Q$ is ASC if and only if the following condition holds.
\begin{enumerate}
\item[(5')] For every $\P\in\Q_{FP}$, $\P\to\Fr$ yields $\P\in\sf SP(\Fr)$.
\end{enumerate}
\end{corollary}

\begin{proof}
Clearly (5') yields condition (5) from Theorem \ref{thm::main characterization}, and hence it implies ASC. For the converse
consider a $\Q$-finitely presented algebra $\P_{\varphi(\bar x)}$ and assume that  it belongs to ${\sf Q}(\Fr)={\sf SPP_U}(\Fr)$. We will show that $\P\in\sf SP(\Fr)$. Strictly, we will prove that for each atomic formula $\psi(\bar x)$ such that $\P_{\varphi(\bar x)}\not\models \psi(\bar v)$  there is a homomorphism $f\colon\P_{\varphi(\bar x)}\to\Fr$ such that $\Fr\not\models\psi(f(\bar v))$. By what we assumed, there is a homomorphism $h\colon\P_{\varphi(\bar x)}\to\Fr^I/U$, for some ultrafilter $U$ over some set $I$, such that $\Fr^I/U\not\models\psi(h(\bar v))$. This means that
\[
\Fr^I/U\models(\exists \bar x)[\varphi(\bar x)\meet\neg\psi(\bar x)]
\]
and, by the elementary equivalence of $\Fr$ with $\Fr^I/U$, there is a tuple of terms $\bar t$ such that
$\Fr\models\varphi(\bar t)\meet\neg\psi(\bar t)$.
Thus we may take as a desired homomorphism $f$ one for which $f(\bar x)=\bar t$ holds.
\end{proof}

From Condition (4) in Theorem \ref{thm::main characterization} we can deduce a supposition under which ASC is equivalent to SC.

\begin{corollary}\label{cor:: ASC <-> SC}
Suppose that every nontrivial algebra from $\Q$ admits a homomorphism into $\Fr$. Then $\Q$ is ASC if and only if it is SC.
\end{corollary}

Note that the assumption of Corollary \ref{cor:: ASC <-> SC} holds when $\Fr$ has an idempotent element, i.e., one element subalgebra. It includes cases of groups or lattices. But in quasivarieties which provide algebraic semantics for  particular deductive systems we rarely have an idempotent element. This is due to the fact that for most encountered cases we have formulas for {\it verum} and {\it falsum} which correspond to two distinct constants in free algebras.
However, even then  Corollary \ref{cor:: ASC <-> SC} is sometimes applicable. It holds e.g. for quasivarieties of  Heyting algebras (Fact \ref{fact::Hey->2}) and  McKinsey algebras, (Lemma \ref{lem:: homo McKinsey to 2}). Note that the latter includes quasivarieties of Grzegorczyk algebras.  We will return to the problem when ASC is equivalent to SC in Proposition \ref{prop:: SC <-> McK} in the case of varieties of closure algebras.

\section{ASC core}\label{sec: asc core}
Let $\Q$ be a variety and $\Fr$ be its free algebra of denumerable rank. Let us consider the interval $[\sf Q(\Fr),\Q]$ in the lattice of subquasivarieties of $\Q$. Notice that all quasivarieties from this interval have the same free algebras. We define the \emph{ASC core} of $\Q$ to be the quasivariety defined relative to $\Q$ by all $\Q$-active quasi-identities and denote it by $\sf ASCC(\Q)$. It follows from the definition of ASC that $\sf ASCC(\Q)$ is the largest ASC quasivariety in $[\sf Q(\Fr),\Q]$. Note that  there does not have to exist a larges ASC subquasivariety of $\Q$, see Example \ref{exm::monounary algebras}.

Since $\sf ASC (\Q)$ is defined relative to $\Q$ by $\Q$-active quasi-identities, $\sf Q(\Fr)$ is defined relative to $\sf ASC (\Q)$ by $\Q$-passive quasi-identities. This fact has a logical interpretation. Namely if a deductive system $\mathcal S$ is ASC, then as a basis of its admissible rules relative to $\mathcal S$ we may take the set of $\mathcal S$-passive rules.

Let us note that $\sf ASC(\Q)$ may be defined also semantically.

\begin{proposition}\label{prop:: ASC core}
For every subalgebra $\C$ of $\Fr$ we have
\[
{\sf ASCC}(\Q)=\{\A\in\Q\mid \A\times\C\in\sf Q(\Fr)\}.
\]
Moreover, a quasivariety $\mathcal R$ from the interval $[\sf Q(\Fr),\Q]$ is ASC if and only if $\mathcal R\leq {\sf ASCC}(\Q)$.
\end{proposition}

\begin{proof}
For the convenience in this proof let us put $\K=\{\A\in\Q\mid \A\times\C\in\sf Q(\Fr)\}$. By Corollary \ref{cor:: about C}, in order to prove that $\K={\sf ASCC}(\Q)$ is is enough to show that $\K\in[\sf Q(\Fr),\Q]$. This means that $\K$ is a quasivariety with $\Fr$ as a free algebra of denumerable rank.
To this end we will check its closeness under $\sf S,P$ and $P_U$ class operators.

So assume first that $\B\leq\A\in\K$. Then $\B\times\C\leq\A\times\C\in\sf Q(\Fr)$, and hence
$\B\in\K$. Now assume that $\A_i\in \K$, for $i\in I$. Then, since
$\sf Q(\Fr)$ is closed under taking direct product, $\left(\prod_{i\in I}\A_i\right)\times \C^I\cong \prod_{i\in I}(\A_i\times \C)\in{\sf Q}(\Fr)$. Since  $\left(\prod_{i\in I}\A_i\right)\times \C$ embeds into $\left(\prod_{i\in I}\A_i\right)\times \C^I$, the algebra $\left(\prod_{i\in I}\A_i\right)\times \C$ also belongs to  $\sf Q(\Fr)$. This proves that $\prod_{i\in I}\A_i\in\K$. For ultraproducts we argue similarly.
Consider an ultrafilter $U$ on a set $I$. Then $\left(\prod_{i\in I}\A_i\right/U)\times \C$ embeds into $\left(\prod_{i\in I}\A_i/U\right)\times \C^I/U\cong \prod_{i\in I}(\A_i\times \C)/U\in{\sf Q}(\Fr)$, and $\prod_{i\in I}\A_i/U\in\K$. In this way we proved that $\K$ is a quasivariety. 

Moreover, the containment $\Fr\times\Fr\in{\sf Q}(\Fr)$ shows that $\K$ has $\Fr$ as a free algebra of denumerable rank.

Now the second statement of the proposition follows from the definition of ASC core or from Corollary \ref{cor:: about C}.
\end{proof}

\section{Projective unification and discriminator varieties}\label{sec:: projective unification}

ASC for varieties (or logics) which are not SC in many  cases was established by means of projective unification, see our historical notes in introduction. Let us look a bit closer at this property (for more details see for instance \cite{Ghi97}). For an equational theory $E$, an $E$-unifier for a finite set $S(\bar x)$ of equations is a substitution $u$, i.e., an endomorphism of a term algebra, such that $(\forall \bar x)\, u(s(\bar x))\=u(t(\bar x))$ belongs to $E$ for every equation $s(\bar x)\=t(\bar x)$ form $S(\bar x)$. However, for our needs it will be more convenient to employ S. Ghilardi algebraic approach \cite{Ghi97}. Let $\V$ be the variety defined by $E$. Instead of working with a finite set of equations $S(\bar x)$, we will deal with the $\V$-finitely presented algebra $\P_{\bigmeet S(\bar x)}$. Then an unifier of $S(\bar x)$ may be identified with an unifier of $\P_{\bigmeet S(\bar x)}$ defined as in Section \ref{sec::quasivariety theory}, i.e., as a homomorphism from $\P_{\bigmeet S(\bar x)}$ into $\Fr$.
A variety $\V$ has \emph{projective unification} if every $\V$-finitely presented unifiable algebra $\P$ is $\V$-projective.
In algebraic terms it means that $\P$ is a retract of $\Fr$. In particular, $\P$ is a subalgebra of $\Fr$. A variety has \emph{unitary unification} if for every $\V$-finitely presented unifiable algebra $\P$ there exists a \emph{most general unifier}, i.e., an unifier through which every unifier of $\P$ can be factorized. Obviously,
projective unification implies unitary unification. 
Note that projective unifiers proved to be very useful in unification and admissibility of rules \cite{Ghi97,Ghi99,Ghi00}.

\begin{corollary}[\protect{\cite{Dzi11}}]\label{cor:: PU->ASC}
If $\sf V(\Q)$ has projective unification, then $\Q$ is ASC.
\end{corollary}

\begin{proof}
It follows directly from Theorem \ref{thm::main characterization} point (4) that $\sf V(\Q)$ is ASC. Now Proposition \ref{prop:: ASC core} yields that $\Q$ is ASC.
\end{proof}

Demonstrating of having projective unification in general has a syntactical nature. However, having projective unification is a stronger property than ASC (see Theorem \ref{thm:: modal companions of Medvedev and Levin} ), and it is not surprising that sometimes it may be established easier, with the aid of semantical methods. We demonstrate this in Example \ref{exp:: discriminator}.

Corollary \ref{cor:: about C} yields that if $\V$ is ASC and $\C\leq\Fr$, then every algebra of the form $\A\times\C$ belongs to $\sf Q(\Fr)$, where $\A\in\V$. On the other hand, if $\V$ has projective unification, then every nontrivial $\V$-finitely presented algebra $\P$ from $\sf Q(\Fr)$ is of the form $\B\times\C$, where $\C\leq\Fr$, in a superficial way, i.e., with $\B$ trivial and $\C\cong\P$. Suppose that $\V$ has projective unification and $\Fr$ has a minimal subalgebra $\C$. Is it then true that every nontrivial (finitely generated or $\V$-finitely presented or just finite) algebra $\A$ in $\sf Q(\Fr)$ have $\C$ as a direct factor? In general: no. We demonstrate this in Example \ref{exm:: S4.3}.
Still, we have the following fact (see Example \ref{exp:: discriminator} for the definition of discriminator variety).

\begin{proposition}[\protect{\cite[Corollary 2.2]{AJN91}}]\label{prop:: simple direct factor}
Suppose $\V$ is a discriminator variety in a finite language. If $\C$ is a finite homomorphic image of a finitely generated member $\A$ of $\V$, then
$\C$ is a direct factor of $\A$.
\end{proposition}

\begin{corollary}\label{cor:: Q(V) when V is discriminator}
Let $\V$ be a discriminator variety in a finite language. Assume that there is a minimal finite subalgebra $\C$ of $\Fr$. Then for every nontrivial $\V$-finitely presented algebra $\P$ the following equivalence holds: $\P\in\sf Q(\Fr)$ if and only if $\P\cong \B\times \C$ for some $\B\in \V$.
\end{corollary}

\begin{proof}
By Corollary \ref{cor:: PU->ASC} and  \cite[Theorem 3.1]{Bur92}, the backward implication holds. 

For the verification take nontrivial $\V$-finitely presented algebra $\P_{\varphi(\bar x)}$ from $\sf Q(\Fr)$. By \cite[Theorems II.9.5]{BS81}, $\P_{\varphi(\bar x)}\in {\sf SPP_U}(\Fr)$ this yields that for every pair of distinct elements from $P_{\varphi(\bar x)}$ there is a homomorphism from $\P_{\varphi(\bar x)}$ into some elementary extension of $\Fr$ separating them. Since $\P_{\varphi(\bar x)}$ is nontrivial, there exists at least one such homomorphism. Note that the algebra from $\V$ admits a homomorphism from $\P_{\varphi(\bar x)}$ iff it satisfies the sentence $(\exists \bar x)\varphi(\bar x)$. Thus there is a homomorphism $h\colon\P_{\varphi(\bar x)}\to\Fr$.   

Clearly, $\Fr$ admits a homomorphism $g$ onto $\C$. Since $\C$ does not have proper subalgebra, $g\circ h$ maps $P_{\varphi(\bar x)}$ onto $C$. Thus, by Proposition \ref{prop:: simple direct factor}, $\C$ is a direct factor of $\P_{\varphi(\bar x)}$.
\end{proof}

\section{Striving for finiteness}\label{sec: criterion}

In order to check the conditions from Theorem \ref{thm::main characterization} it is possible that one has to work on infinite algebras. The following question arises: Under what conditions can we simplify verification of ASC by restricting condition (3) from Theorem \ref{thm::main characterization} to finite algebras? In this section we will propose a solution to this problem, namely Theorem \ref{thm::criterion}. In the next section we will show some of its applications.

Let us start with recalling needed notions.
We say that a class $\mathcal K$ of algebras has \emph{finite model property} (FMP for short) if ${\sf V}(\mathcal K)$ is generated, as a variety, by finite members from $\mathcal K$. Note that it may happen that a quasivariety does not have FMP while the variety it generates does. A class $\mathcal K$ has \emph{strong finite model property} (SFMP for short) if $\sf Q(\mathcal K)$ is generated, as a quasivariety, by finite members from $\mathcal K$. In particular, every locally finite (with all finitely generated algebras being finite) quasivariety has SFMP. A quasivariety $\Q$ has  \emph{equationally definable principal relative congruences}  (EDPRC for short and EDPC for varieties) if there is a finite family of equations $s_k(u,v,x,y)\=t_k(u,v,x,y)$, $k\leq n$, such that for every $a,b,c,d\in A$ and $\A\in \Q$
\[
(c,d)\in\theta_\Q(a,b)\quad\text{iff}\quad \A\models \bigmeet_{k\leq n} s_k(c,d,a,b)\=t_k(c,d,a,b).
\]

\begin{theorem}\label{thm::criterion}
Let $\Q$ be a quasivariety in a finite language with FMP and EDPRC. Assume that $\Fr$ has a finite $\Q$-simple subalgebra $\C$. Then $\Q$ is ASC if and only if for every finite $\Q$-subdirectly irreducible algebra $\S$ we have
\[
\S\leq \Fr\quad\text{ or }\quad\S\times\C\leq \Fr.
\]
\end{theorem}

Let us emphasize that the assumptions of Theorem \ref{thm::criterion} are very natural from the perspective of logic. Indeed, assume that $\Q$ gives an algebraic semantics for a deductive system $\mathcal S$. Then having FMP by $\Q$ with a recursively enumerable axiomatization of $\sf V(\mathcal K)$ yield the decidability of the equational theory of $\mathcal K$ and hence the decidability of $\Th(\mathcal S)$ \cite[Theorem 3]{McK43}. Furthermore, having EDPRC by $\Q$ corresponds to deduction-detachment theorem for $\mathcal S$ \cite[Theorem 5.5]{BP01}, \cite[Theorem 4.6.13]{Cze01}. The algebra $\C$ may be often chosen as an algebra with elements which correspond to {\it verum} and {\it falsum}.

\begin{lemma}\label{lem:: SFMP + finite SI ok -> ASC}
Assume that $\Q$ has SFMP and $\C$ is a subalgebra of $\Fr$. If for every finite $\S\in\Q_{SI}$, $\S\times\C\in\sf Q(\Fr)$, then $\Q$ is ASC.
\end{lemma}

\begin{proof}
By Proposition \ref{prop:: ASC core}, the class of finite $\Q$-subdirectly irreducible algebras is contained in ${\sf ASCC}(\Q)$. Consequently, by Proposition \ref{prop::Mal'cev thm}, all finite algebras from $\Q$ are in ${\sf ASCC}(\Q)$. Thus, by SFMP, $\Q={\sf ASCC}(\Q)$. This means that $\Q$ is ASC.
\end{proof}

For a congruence $\alpha$ of $\A$ and $\beta$ of $\B$ let $\alpha\times\beta$ be a congruence of $\A\times\B$ given by $\{((a_1,b_1),(a_2,b_2))\in (A\times B)^2\mid (a_1,b_1)\in\alpha\text{ and }(a_2,b_2)\in\beta\}$. A quasivariety $\Q$ has \emph{Fraser-Horn property} (FHP for short) if for every algebras $\A,\B$ each $\Q$-congruence of the product $\A\times\B$ decomposes as $\alpha\times\beta$, where $\alpha$ is a $\Q$-congruence of $\A$ and $\beta$ is a $\Q$-congruence of $\B$.  Every relative congruence distributive quasivariety has FHP but this notion is more general, see \cite{CzD92}.

\begin{lemma}\label{lem:: ASC-> finite SI are ok}
Assume that $\Q$ is a quasivariety in a finite language
which has FHP and $\Fr$ has a finite $\Q$-simple subalgebra $\C$. If $\Q$ is ASC then for every finite $\S\in\Q_{SI}$ we have
\[
\S\leq \Fr\quad\text{ or }\quad\S\times\C\leq \Fr.
\]
\end{lemma}

\begin{proof}
Let $\S$ be a finite $\Q$-subdirectly irreducible algebra. By ASC, $\S\times\C\in\sf Q(\Fr)={\sf SPP_U}(\Fr)$. This means that for each pair of distinct elements in $S\times C$ there is a homomorphism from $\S\times\C$ into an ultrapower of $\Fr$ that separates them. Let $(a,b)\in S^2$ be a pair which belongs to every $\Q$-congruence of $\S$ that is not the equality relation $\id_S$ on $S$. Further, let $c$ be an element of $\C$. Let $h\colon \S\times\C\to \G$ be a homomorphism such that $h(a,c)\neq h(b,c)$, where $\G$ is an ultrapower of $\Fr$. Then FHP yields that $\ker(h)=\alpha\times\beta$, where $\alpha$ is a $\Q$-congruence of $\S$ and $\beta$ is a $\Q$-congruence of $\C$. As $h(a,c)\neq h(b,c)$, $\alpha$ equals $\id_S$ and, by $\Q$-simplicity of $\C$, $\beta$ equals $\id_C$ or $C^2$. Thus, either $\S$ or $\S\times\C$ embeds into $\G$. By the finiteness of the language of $\Q$ and the finiteness of both algebras, at least one of them embeds \mbox{into $\Fr$}.
\end{proof}

We need two facts from the literature.

\begin{proposition}[\protect{\cite[Theorem 3.3]{BvA02}}]\label{prop::Blok and van Alten}
For a quasivariety FMP and EDPRC yields SFMP.
\end{proposition}

\begin{proposition}[\protect{\cite[Theorem 5.5]{BP01}, \cite[Theorem Q.9.3]{Cze01}}]\label{prop:: EDPRC->FHP}
A quasivariety with EDPRC is relative congruence-distributive, and thus has FHP.
\end{proposition}

\begin{proof}[Proof of Theorem \ref{thm::criterion}]
For the backward direction combine Proposition \ref{prop::Blok and van Alten} and Lemma \ref{lem:: SFMP + finite SI ok -> ASC}. For the forward direction combine Proposition \ref{prop:: EDPRC->FHP} and Lemma \ref{lem:: ASC-> finite SI are ok}.
\end{proof}

As a matter of fact, there is an analog of Theorem \ref{thm::criterion} for SC.

\begin{corollary}\label{cor:: skonczona characterizacja SC}
Let $\Q$ be a quasivariety in a finite language with EDPRC. Then $\Q$ is SC if and only if every finite $\Q$-subdirectly irreducible algebra is a subalgebra of $\Fr$.
\end{corollary}

\begin{proof}
The backward direction follows from Proposition \ref{prop::Blok and van Alten} and the fact that  all finite algebras from $\Q$ are in $\sf Q(\Fr)$. This fact follows form Proposition \ref{prop::Mal'cev thm} and the assumption. The forward implication may be proved similarly, but easier, as Lemma \ref{lem:: ASC-> finite SI are ok}.
\end{proof}

\begin{remark}
Corollary \ref{cor:: skonczona characterizacja SC} was obtained in \cite[Theorem 5.1.8]{Ryb97} under some additional condition. But in the cases of intermediate logics and of normal extensions of K4 modal logics \cite[Corollary 5.1.10]{Ryb97} the formulation presented there is the same as ours.
\end{remark}

Several forms of definability of relative principal congruences which are weakenings of EDPRC were considered in the literature. They correspond to variants of deduction-detachment theorem for deductive systems. Among them the property of having equationally semi-definable principal relative congruences, corresponding to contextual deduction-detachment theorem \cite[Theorem 9.2]{Raf11}, proves to be sufficient for Theorem \ref{thm::criterion} to work. Indeed, having equationally semi-definable principal relative congruences yields relative congruence-distributivity, and with FMP yields SFMP \cite[Theorem, 8.7, Corollary 3.7]{Raf11}

\begin{problem}
Is it possible to weaken the assumption of Theorem \ref{thm::criterion} of having EDRPC to, having relative congruence extension property, corresponding to local deduction-detachment theorem \cite[Corollary 3.7]{BP88}, or to having parameterized equationally definable principal relative congruences, corresponding to parameterized deduction-detachment theorem \cite[Section 2.4]{Cze01}?
\end{problem}

\section{Examples}\label{sec:: examples}

In this section we will give several examples of ASC varieties. The main objective is to present varieties which characterize ASC$\setminus$SC logics. The exception is given by varieties of monounary algebras and varieties of
bounded lattices. They are intended to illustrate how one may apply Theorems \ref{thm::main characterization} and \ref{thm::criterion} and the techniques used in their proofs. Also the example of monounary algebras shows that there does not have to to exists a largest ASC subquasivariety of a given quasivariety. Moreover, the example of bounded lattices shows that some ``plausible'' condition for ASC is actually strictly weaker than ASC.

We will use a nonstandard notation for operations in algebras and instead of $\join,\meet,\to,\neg$ we will write $\mjoin,\mmeet,\Rightarrow,\mneg$ symbols. We do so in order to make a visible distinction between a language and the meta-language.

\begin{example}{\bf Monounary agebras.}\label{exm::monounary algebras}
Let $\V$ be the class of all monounary algebras. These are algebras with just one basic operation, denoted by $f$, which is unary. We claim that ${\sf ASCC}(\V)$ is defined by the quasi-identity
\[
j=(\forall x,y)[f(x)\=f(y)\to x\=y].
\]
We may identify a free monounary algebra $\Fr(1)$  with $(\Nat,f\colon x\mapsto x+1)$. Note that $\Fr$ is isomorphic with a disjoint union of denumerable many copies of $\Fr(1)$. We clearly have $\Fr\models j$ and $\Fr\not\models j^*$. Hence ${\sf ASCC}(\V)\models j$. Now in order to prove our claim it is enough to show that the quasivariety defined by $j$ is ASC. To this end one may use the condition (2') from Corollary \ref{cor:: about C}. Indeed, if $\A\models j$, then $\A\times\Fr(1)$ is a disjoint union of subalgebras generated by $(a,n)$, $a\in\A$, $n\in\Nat$, where $a\not\in f(A)$ or $n=0$. Each of these subalgebras is isomorphic to $\Fr(1)$. Therefore  $\A\times\Fr(1)$ is free for $\V$ and belongs to $\sf Q(\Fr)$ (actually, all nontrivial members of $\sf Q(\Fr)$ are free for $\V$).

Now consider a variety $\W$ defined by $(\forall x,y)[f(x)\=f(y)]$. Then it has, up to isomorphism, only one subdirectly irreducible algebra $(\{0,1\},x\mapsto \min(1,x+1))$. Moreover,  this algebra embeds into every nontrivial member of $\W$. Thus $\W$ is a minimal quasivariety and is SC. But $j$ is not valid in $\W$ and $\W\not\subseteq{\sf ASCC}(\V)$. This shows that there does not have to exist a largest (A)SC subquasivariety of a given quasivariety.
\qed
\end{example}

\begin{example}{\bf Varieties of bounded lattices.}

By a \emph{bounded lattice} we mean an algebra $\mathbf L$ with a lattice reduct and with two constants $0$ and $1$ which are the bottom and the top elements in $\mathbf L$ respectively.

Due to the lack of FMP, Theorem \ref{thm::criterion} does not apply to all varieties of bounded lattices. (For instance the variety defined by modularity law does not have FMP. In \cite{Fre79} an identity $e$ was found that holds in all finite modular lattices, but does not hold is some infinite one $\bf L$.  Clearly, $e$ also holds in all finite bounded lattices, and does not hold in the bounded expansion of $\bf L$.)
Still, the argument from the  proof may be used to show that there are only two ASC (SC in fact) varieties of bounded lattices. Strictly, the proof of Lemma \ref{lem:: ASC-> finite SI are ok} yields also the following fact.

\begin{lemma}\label{lem:: ASC + distributivity -> prawie ok}
Assume that $\Q$ is a quasivariety with FHP and that $\Fr$ has a finite $\Q$-simple subalgebra $\C$. If $\Q$ is ASC then for every  $\S\in\Q_{SI}$ we have
\[
\S\leq \G\quad\text{ or }\quad\S\times\C\leq \G.
\]
for some ultrapower $\G$ of $\Fr$.
\end{lemma}

Here by $\2$ we denote the bounded lattice $(\{0,1\},\mmeet,\mjoin,0,1)$. Note that $\2$ is free of rank zero for every nontrivial variety of bounded lattices. The following lemmas are folklore.

\begin{lemma}\label{lem:: SI bounded lattices}
Let $\S$ be a subdirectly irreducible bounded lattice not isomorphic to $\2$. Then $1$ does not have the unique lower cover in $\S$. In particular, if $\S$ is finite,  $1$ is join-reducible in $\S$.
\end{lemma}

\begin{proof}
Assume  that there is the unique lower cover $1_*$ of $1$ in $\S$. Consider two congruences of $\S$
\begin{align*}
\alpha &=\{1,1_*\}^2\cup \id_S,\\
\beta &=\{a\in S\mid 0\leq a\leq 1_*\}^2\cup\id_S.
\end{align*}
Then $\alpha>\id_S$ and $\alpha\cap\beta=\id_S$. Thus, by subdirect irreducibility, $\beta=\id_S$ and $\S$ must be isomorphic to $\2$. This leads us to a contradiction with our assumption.
\end{proof}

\begin{lemma}\label{lem:: 1 is join irr in free bounded lat}
Let $\V$ be a nontrivial variety of bounded lattices and $\G$ be an ultrapower of $\Fr$. Then  $1$ is join-irreducible in $\G$.
\end{lemma}

\begin{proof}
Since the conclusion of the lemma is expressible by a first order sentence and $\G$ is  elementarily equivalent to $\Fr$, it is enough to prove it for $\Fr$.

Assume that $1$ is join-reducible in $\Fr$. Then there are
$p(\bar v),q(\bar v)\in F$ such that  $p,q <1$ and $p\mjoin q=1$ in $\Fr$. This yields that $\V$ satisfies the identity $\hat p\mjoin \hat q\=1$, where $\hat p=p(0,\ldots,0),\hat q=q(0,\ldots,0)\in F(0)=\{0,1\}$. In particular, in $\Fr(0)$, which is isomorphic to $\2$, we have
$\hat p\mjoin \hat q=1$. Thus at least one of $\hat p,\hat q$, say $\hat p$, equals 1. Since in bounded lattices all term operations are monotone, $\hat p\leq p$. Hence $p=1$ in $\Fr$. This gives a contradiction.
\end{proof}

Let $\N_5$ be a 5-element lattice in which non-top and non-bottom elements form a disjoint union of an element with a two-element chain (exactly two among these elements are comparable). Let $\M_3$ be a 5-element lattice in which non-top and non-bottom elements form a three-element antichain (all of them are incomparable). Let $\N_5^b$ and $\M_3^b$ be bounded lattices with the lattice reducts $\N_5$ and $\M_3$ respectively. By \emph{distributivity law} we mean the identity
\[
(\forall x,y,z)\;[x\mmeet(y\mjoin z)\=(x\mmeet y)\mjoin(x\mmeet z)].
\]

\begin{lemma}\label{lem:: bounded distributive lattices characterization}
Let $\V$ be a variety of bounded lattices. Then the following conditions are equivalent:
\begin{enumerate}
\item $\V$ satisfies distributivity law,
\item $\V=\sf V(\2)$ or $\V$ is the trivial variety,
\item $\N^b_5\not\in\V$ and $\M^b_3\not\in\V$.
 \end{enumerate}
\end{lemma}

\begin{proof}
\mbox{}\\
\noindent(1)$\Rightarrow(2)$ It follows from Pristley duality \cite[Theorem 11.23]{DP02} that the class of distributive bounded lattices coincide with ${\sf SP}(\2)$. Thus, since $\2$ is distributive, ${\sf SP}(\2)=\sf V(\2)$. On the other hand, every nontrivial bounded lattice has a subalgebra isomorphic to $\2$. Hence there are only two varieties of bounded lattices satisfying distributivity law: $\sf V(\2)$ and the trivial variety.

\noindent (2)$\Rightarrow$(1)$\Rightarrow$(3) It is routine.

\noindent(3)$\Rightarrow(1)$ Assume that in $\V$ there is a non-distributive bounded lattice $\mathbf L$. Then its lattice reduct has a sublattice $\mathbf K$ isomorphic to $\M_3$ or $\N_5$ \cite[Theorem I4.10]{DP02}. Let $\mathbf K^b$ be a bounded sublattice of $\mathbf L$ generated by $K$. Note that the lattice reduct of $\mathbf K^b$ may differ from $\mathbf K$ only by having an additional element on the top and/or having an additional element in the bottom. In either case, $\mathbf K^b$ has one of bounded lattices $\M_3^b$, $\N_5^b$ as a homomorphic image. Thus one of these algebras belongs to $\V$.
\end{proof}

\begin{proposition}\label{prop:: bounded lattices}
Let $\V$ be a variety of bounded lattices. Then the following conditions are equivalent:
\begin{enumerate}
\item $\V$ is SC,
\item $\V$ is ASC,
\item $\V$ satisfies distributivity law.
\end{enumerate}
\end{proposition}

\begin{proof}
\mbox{}\\
\noindent(1)$\Rightarrow(2)$ It is obvious.

\noindent(2)$\Rightarrow(3)$ Assume that in $\V$ distributivity law does not hold. Then, by Lemma \ref{lem:: bounded distributive lattices characterization},  at least one of $\M_3^b$, $\N_5^b$  belongs to $\V$. For convenience, let us denote it by $\S$. Clearly, $\S$ is finite, subdirectly irreducible and has the top element join-reducible. Since in $\S\times\2$ the top element is also join-reducible, Lemma \ref{lem:: 1 is join irr in free bounded lat} yields that neither $\S$ nor $\S\times\2$ embeds into any ultrapower of $\Fr$. Thus, by Lemma \ref{lem:: ASC + distributivity -> prawie ok}, $\V$ cannot be ASC.

\noindent(3)$\Rightarrow(1)$ Lemma \ref{lem:: bounded distributive lattices characterization} tells us that there are only two varieties of distributive bounded lattices: the trivial one, which is clearly SC, and the minimal one $\sf V(\2)$. In fact $\sf V(\2)$ is minimal also as a quasivariety and as such must be also SC.
\end{proof}

We finish this example by one remark.
Consider the condition for quasivarieties obtained by syntactic mixing the conditions from Theorem \ref{thm::main characterization}.
\begin{enumerate}
\item[(6)] For every $\S\in\Q_{SI}$, $\S\to\Fr$ yields $\S\in\sf Q(\Fr)$.
\end{enumerate}
Theorem \ref{thm::main characterization} shows that (6) follows from ASC. But is not equivalent to ASC. In order to see this, let us consider  the variety $\V$ generated by $\M_3^b$. By Proposition \ref{prop:: bounded lattices}, $\V$ is not ASC. Let us check that nevertheless the condition (6) is fulfilled.  There are, up to isomorphism, exactly two subdirectly irreducible algebras in $\V$: $\2$ and $\M_3^b$. Clearly $\2\in{\sf Q}(\Fr)$. Moreover, $\M_3^b$ does not admit a homomorphism into $\Fr$. Indeed, since $\M_3^b$ is simple, a homomorphic image of $\M_3^b$ would have just one element, which is impossible since $\Fr$ does not have idempotents, or be isomorphic to $\M_3^b$, which is also impossible as we showed in the proof of Proposition \ref{prop:: bounded lattices}.

Note that the unbounded case is different. In particular, $\sf V(\M_3)$ is SC \cite{Igo74}.
\qed
\end{example}

Let us move to examples that come from logic. We are mainly interested in normal modal logics, and, in particular, in normal extensions of transitive and reflexive modal logic S4. Every such extension has an adequate semantics given by a variety of closure algebras \cite{MT48}, \cite[Chapter 10]{DH01}. An algebra $\M$ is a \emph{modal algebra} if it has a Boolean algebra reduct and beside Boolean operations one unary operations $\Diamond$ such that for all $a,b\in M$
\[
\Diamond 0=0,\quad \Diamond(a\mjoin b)=\Diamond a\mjoin \Diamond b.
\]
If in addition for every $a\in M$ it satisfies
\[
a\leq \Diamond a=\Diamond\Diamond a
\]
we call it a \emph{closure algebra}.
Let $\mBox x$ $=\mneg\Diamond\mneg x$. Element $a$ of a closure algebra is \emph{closed} (\emph{open}) if $a=\Diamond a$ ($a=\mBox a$ respectively). We picture a closure algebra $\M$ by drawing the Hasse diagram of the ordered set $(M,\leq)$, where $\leq$ is given by the lattice structure of $\M$.  We draw closed elements as $\Diamond$, open as $\mBox$, open and closed as $\DiamondBox$, and others as $\piccirc$.
The simplest nontrivial closure algebra, denoted by $\2$ and depicted in Figure \ref{fig:: 2 and S_2}, has two elements and $\Diamond$ operation acts on it identically. In particular it is term equivalent to a two-element Boolean algebra. It is important to note that $\2$ embeds into every nontrivial closure algebra. Moreover, $\2$ is free of rank zero for every nontrivial variety of closure algebras.

Let us recall that congruences of a closure algebras $\M$ are with one to one correspondence with \emph{open filters}, i.e., Boolean filters which are additionally closed under $\mBox$ operation. Strictly, for a congruence $\alpha$ its corresponding open filter is the class $1/\alpha$. From this one can see that the variety of closure algebras is congruence-distributive. Actually a stronger statement is true: every variety of closure algebras has EDPC witnessed by the equation $\mBox(x\Leftrightarrow y)\mto(u\Leftrightarrow v)\=1$.
Note also that each element $a$ of a closure algebra $\M$ which is open and closed gives a direct product decomposition $\M\cong\M/\alpha\times\M/\beta$, where $\alpha$ is a congruence generated by $(1,a)$ and $\beta$ is a congruence generated by $(1,\mneg a)$.

\begin{example}\label{exp:: monadic algebras}{\bf Varieties of monadic algebras}. A closure algebra $\M$ is a \emph{monadic algebra} if for all $a\in M$ we have
\[
\Diamond \mBox a=\mBox a.
\]
This means that all open elements in $\M$ are also closed.
Recall that varieties of monadic algebras form adequate semantics for normal extensions of transitive, reflexive and symmetric S5 modal logic \cite{MT48}, \cite[Chapter 10]{DH01}. As we already noted in the introduction, every variety of monadic algebras is a discriminator variety. Hence it has projective unification and is ASC. Since the variety of monadic algebras is for us a prototypical example of an ASC variety which is not SC, let us look at it from an algebraic perspective.
For this purpose we will need to recall basic facts about monadic algebras.

For a positive integer $l$ let $\S_l$ be the closure algebra with $l$ atoms and with 0 and 1 as the only closed elements. The algebra $\S_1$, which is isomorphic to $\2$, and the algebra $\S_2$ are depicted in Figure \ref{fig:: 2 and S_2}. Clearly, all $\S_l$ are monadic. Let $\V$ be a variety of monadic algebras.
Then $\V$ is semisimple, i.e., all its subdirectly irreducible algebras are simple. Moreover, every finite simple closure algebra is isomorphic to one of $\S_l$ \cite[Lemma 8, Theorem 7]{Hal56}, \cite[Theorem 4.2]{KQ76}.
This gives that every finite monadic algebra $\M$ is isomorphic to a product of those $\S_l$ which are its homomorphic images. Indeed, every maximal congruence $\alpha$ of a monadic algebra $\M$ is generated by a pair $(a,1)$ where $a$ is open, and hence also closed, in $\M$. Thus $\M$ is isomorphic to the product $\M/\alpha\times \M/\beta$, where $\beta$ is generated by $(\mneg a,1)$. Since $\V$ is locally finite \cite{Bas58}, this applies to $\Fr(k)$ for every \mbox{finite  $k$}. So we have
\[
\Fr(k)\cong\prod_{l=1}^{m}\S_l^{d_l}
\]
for some natural numbers $m,d_1,\ldots,d_m$. Note also that if $\S_l\in\V$ and $k\geq l$, then $\S_l$ is a homomorphic image of $\Fr(k)$ and $d_l\geq 1$.
(An exact structure of $\Fr_\W(k)$ may be deduced from \cite{Bas58,Hal59,KQ76} where free monadic algebras are described.)

Let us use Theorem \ref{thm::criterion} in order to show that $\V$ is ASC. As we already noted $\V$ has EDPC and, since it is locally finite, it has FMP. Moreover a two-element closure algebra $\2$ embeds into every nontrivial monadic algebra. Thus the assumptions of Theorem \ref{thm::criterion} hold. Let us verify the condition from the theorem. For a trivial $\V$ it vacuously holds. So assume that $\S_l\in\V$ for some positive integer $l$. Take $n\geq l$. Then, according to what we already wrote, $\Fr(n)\cong \S_l\times\M$ for some monadic algebra $\M$, and hence $\S_l\times\2$ embeds into $\Fr(n)$ (when $\M$ is nontrivial) or $\S_l$ embeds into $\Fr(n)$ (when $\M$ is trivial).

Note that there are only two SC varieties of monadic algebras, namely the trivial one and $\sf V(\2)={\sf SP}(\2)$. The latter one is actually term equivalent to the variety of Boolean algebras. Indeed, all other varieties of monadic algebras contain $\S_2$. Thus, as we indicated in the introduction and will prove in Proposition \ref{prop:: SC <-> McK}, they cannot be SC.
\qed\end{example}

\begin{example}\label{exp:: discriminator}{\bf Locally finite discriminator varieties.}
Actually, the argument for ASC from the previous example may be used in a more general setting. Recall that a variety $\V$ is a  \emph{discriminator variety} if it is generated by a class $\mathcal K$ of algebras for which there is a term $t(x,y,z)$ such that for all $a,b,c\in A$, $\A\in\mathcal K$ we have
\[
t(a,b,c)=
\begin{cases}
a & \text{ if } a\neq b\\
c & \text{ if } a= b
\end{cases}.
\]
Assume that $\V$ is a locally finite discriminator variety in a finite language and that there exists an algebra $\C$ that embeds into every nontrivial member of $\V$. Then $\V$ has FMP and EDPC \cite[Page 200]{BP82}. These assumptions are met in e.g. in the varieties of monadic algebras \cite{KQ76}, of MV${}_n$-algebras \cite{CDM00}, in locally finite varieties of relation algebras \cite{Mad06} or diagonal free cylindric algebras \cite{HMT85}. By Proposition \ref{prop:: simple direct factor}, every subdirectly irreducible (which is actually here the same as simple) algebra $\S$ in $\V$ is a direct factor of $\Fr(n)$ for $n\geq|S|$. Thus if $\S$ is not isomorphic to $\Fr(n)$, then it is a proper direct factor of $\Fr(n)$ and then $\S\times\C$ embeds into  $\Fr(n)$. Thus the assumptions  and the condition form Theorem \ref{thm::criterion} hold. Therefore $\V$ is ASC.
\qed\end{example}

Let us move to a more complicated examples, varieties of closure algebras which are not discriminator.

\begin{example}\label{exm:: S4.3}{\bf Varieties of S4.3-algebras.}
Let $\V_{\rm S4.3}$ be the variety generated by the closure algebras in which open elements form a chain. Alternatively one may define $\V_{\rm S4.3}$, relative to the variety of closure algebras, by
\[
(\forall x,y)[\mBox(\mBox x\mto y)\mjoin\mBox(\mBox y\mto x)].
\]
Note that  $\V_{\rm S4.3}$ characterizes the modal logic $S4.3$, see e.g. \cite{CZ97}. Let $\V$ be a subvariety of $\V_{\rm S4.3}$. We already noted in the introduction that $\V$ is ASC. Let us now argue for it without projective unification. By Bull theorem \cite{Bul66}, $\V$ has FMP. Thus the assumptions of Theorem \ref{thm::criterion} are satisfied for $\V$. Moreover, the condition from the theorem is verified in \cite[Lemma 2]{Ryb84}. Thus $\V$ is ASC. Note that Rybakov in \cite[Theorem 5]{Ryb84} obtained the quasi-equational base for $\sf Q(\Fr)$.

In Section \ref{sec:: projective unification} we formulated the problem whether all finite/finitely presented/finitely generated nontrivial algebras in $\sf Q(\Fr)$ have a minimal subalgebra of $\Fr$ as a direct factor provided $\V$ is a variety with projective unification. To falsify this let us consider the closure algebra $\M$ depicted in Figure \ref{fig:: some S4.3 algebra} and the subvariety $\V$ of $\V_{\rm S4.3}$ containing $\M$. \begin{figure}[h]
\[
\xy
(10,30)*{\DiamondBox}="ONE";
(0.01,19.81)*{\Diamond}="AB"; (10,20)*{\piccirc}="AC"; (20,20)*{\piccirc}="BC";
(0,10)*{\piccirc}="A"; (10,10)*{\piccirc}="B"; (20,10)*{\mBox}="C";
(10,0)*{\DiamondBox}="ZERO";
(9.07,1.07); (.7,9.44) **\dir{-}; 
(10,1.1); (10,9.2) **\dir{-}; 
(10.93,1.07);  (19.1,9.24) **\dir{-}; 
"A"; (0,19.0) **\dir{-}; 
(9.37,19.52); (.60,10.75) **\dir{-}; 
(9.37,10.78); (0.60,19.55) **\dir{-}; 
(19.37,19.52); (10.60,10.75) **\dir{-}; 
(19.07,11.07); (10.7,19.44) **\dir{-};  
(20,11.1); (20,19.2) **\dir{-}; 
(9.07,29.22); (.60,20.75) **\dir{-}; 
"AC"; (10,29.2) **\dir{-}; 
(19.37,20.78); (10.95,29.20) **\dir{-}; 
\endxy
\]
\caption{The algebra $\M$ from $\V_{\rm S4.3}$.}
\label{fig:: some S4.3 algebra}
\end{figure}
Then $\M$ has $\2$ as a homomorphic image. Thus, by Theorem \ref{thm::main characterization} point (4), $\M\in \sf Q(\Fr)$. Nevertheless, $\2$ is not a direct factor of $\M$.
\qed\end{example}

Now we will move to varieties of Heyting algebras. A \emph{Heyting algebra} (called sometimes a pseudo-Boolean algebra) $\H$ is a bounded lattice expanded by one binary operation $\mto$ such that for all $a,b,c\in H$
\[
a\mmeet c\leq b \quad\text{iff}\quad c\leq a\mto b.
\]
Let $\mneg x=x\mto 0$.
Varieties of Heyting algebras constitute an adequate semantics for intermediate logics. In particular, the class of all Heyting algebras, which turns out to be a variety, characterizes intuitionistic logic \cite[Chapter 7]{CZ97}. As in the case of closure algebras, there is exactly one minimal (quasi)variety of Heyting algebras. It is generated by a two-element Heyting algebras, again denoted as $\2$.

Corollary \ref{cor:: ASC <-> SC} and Fact \ref{fact::Hey->2} yield that for varieties of Heyting algebras SC is in fact equivalent to ASC. Nevertheless, they are strongly connected to varieties of closure algebras. In the next section we will show how to construct, starting from an SC variety of Heyting Algebras, infinitely many varieties of closure algebras which are ASC$\setminus$SC.

\begin{example}\label{exp:: Levin and Medvedev}{\bf Levin and Medvedev varieties.}
Recall that with every ordered set $\mathbf O$,  the Heyting algebra $\mathbf O^+$ of its up-directed subsets is associated. Then $\mathbf O$, treated as an intuitionistic frame, validates the intuitionistic formula $t(\bar x)$ iff $\mathbf O^+$ validates the identity $(\forall\bar x)[t(\bar x)\=1]$, see e.g. \cite[Chapter 7]{CZ97}. For a natural number $n$ let $(2,\leq)^n$ be the power of the ordered set $(2,\leq)$ with $2=\{0,1\}$ and $0\leq 1$. Let $\Lev_n$ be the ordered set obtained from $(2,\leq)^n$ by removing the top element. Since the algebra $\Lev_2^+$ will be important in our investigations, we will use a more intuitive notation $\2^2\oplus\1$ for it.
Note that the logic characterized by all $\Lev_n$ is Medvedev finite problems logic \cite{Med62,med66}, and an intermediate logic is one of Levin logics \cite{Lev69} iff it is Medvedev logic or it is characterized by one of frames $\Lev_n$ \mbox{\cite[Theorem 3.1]{Skv99}}.
\begin{figure}[h]
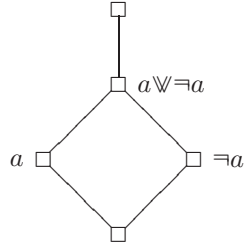

\[
\xy
(10,30)*{\mBox}="ONE";
(10,20)*{\mBox}="AC";
(0,10)*{\mBox}="A";  (20,10)*{\mBox}="C";
(10,0)*{\mBox}="ZERO";
(-3.5,10)*{a}; (24.6,10)*{\mneg a}; (17,20)*{a\mjoin\mneg a};
(9.07,1.07); (.9,9.24) **\dir{-}; 
(10.93,1.07);  (19.1,9.24) **\dir{-}; 
(9.07,19.22); (.90,11.05) **\dir{-}; 
(19.07,11.07); (10.9,19.24) **\dir{-};  
(10,21.1); (10,29.2) **\dir{-}; 
\endxy
\]
\caption{The Heyting algebra $\2^2\oplus\1$.}
\label{fig:: Lev_2^+}
\end{figure}

For $n\in\Nat$ let
\[
\V_{Lev_n}={\sf V}(\Lev_n^+)
\]
and $\V_{Med}$ be their varietal join
\[
\V_{Med}=\bigjoin_{n\in\Nat}\V_{Lev_n}.
\]

The following basic property of Heyting algebras, which may be deduced from e.g. \cite[Statement VI.6.5]{Ras74}, will be needed.

\begin{fact}\label{fact::Hey->2}
Let $\H$ be a Heyting algebra and $a$ its non-zero element. Then there is homomorphism $h\colon \H\to \2$ such that $h(a)=1$.
\end{fact}

\begin{lemma}\label{lem:: 2^2 finitely presented}
Let $\V$ be a nontrivial variety of Heyting algebras. Then $\2^2$ is a $\V$-finitely presented unifiable algebra.
\end{lemma}

\begin{proof}
Since the language of $\V$ is finite and $\2^2$ is a finite algebra belonging to $\V$, Fact \ref{fact::finite->fp} yields that $\2^2$ is $\V$-finitely presented. One can also show that $\2^2$ is isomorphic to $\P_{x\mjoin \mneg x\= 1}$.

The unifiability follows from the existence of a projection from $\2^2$ onto $\2$. Note that $\2$ is free for every nontrivial variety of bounded lattices. 
\end{proof}

\begin{lemma}\label{lem:: 2^2 not embeddable}
Let $\V$ be a variety of Heyting algebras containing $\2^2\oplus\1$. Then $\2^2$ does  not embed into $\Fr$.
\end{lemma}

\begin{proof}
Striving for a contradiction, suppose that $\2^2$ embeds into $\Fr$. Then there is $t\in F$ such that
\[
0<t<1,\; 0<\mneg t<1,\; t\mjoin \mneg t=1,\; t\mmeet \mneg t=0.
\]
Note that neither $t$ nor $\mneg t$ is a Boolean tautology. Indeed, by Fact \ref{fact::Hey->2}, there is a homomorphism $k\colon\Fr\to \2^2$ such that $k(t)=(1,0)$ and $k(\mneg t)=(0,1)$. The Heyting algebra $\2^2\oplus\1$ has  two atoms $a=\{(1,0)\}$, $\mneg a=\{(0,1)\}$, and one coatom $a\mjoin\mneg a=\{(0,1),(1,0)\}$. Let $g\colon V\to 2^2\oplus 1$ be a mapping given by
\[
g(v)=
\begin{cases}
a &\text{ if } k(v)=(1,0)\\
\mneg a &\text{ if } k(v)=(0,1)\\
a\mjoin\mneg a&\text{ if } k(v)=(1,1)\\
0 &\text{ if } k(v)=(0,0)
\end{cases},
\]
and  $\bar g\colon\Fr\to\2^2\oplus\1$ be the homomorphic extension of $g$. Let $h\colon\2^2\oplus\1\to\2^2$ be a surjective homomorphism that maps $a\mjoin\mneg a$ onto 1 and $a$ onto $(1,0)$. Note that $h^{-1}(1)=\{1,a\mjoin\mneg a\}$ is the only its coset containing more than 1 element. We have
$k|_V=h\circ g$ and hence, by the universal mapping property of $\Fr$, $k=h\circ\bar g$. Therefore $\bar g(t)=a$ and $\bar g(\mneg t)=\mneg a$. Now we compute in $\2^2\oplus\1$
\[
1=\bar g(1)=\bar g(t\mjoin\mneg t)=\bar g(t)\mjoin\mneg \bar g(t)=a\mjoin\mneg a<1.
\]
This leads to a contradiction.
\end{proof}

Note that if $\V$ is a variety of Heyting algebras containing a three-element algebra, then $\2^2$ is not $\V$-projective, see the last remark in \cite{Ghi99}.

\begin{lemma}\label{lem:: no unitary unification}
Let $\V$ be an ASC variety and  $\P\in\V$ be a $\V$-finitely presented unifiable algebra that does not embed into $\Fr$. Then $\V$ cannot have unitary unification.
\end{lemma}

\begin{proof}
By Corollary \ref{cor:: on finitely presented}, $\P\in\sf SP(\Fr)$. This means that there are unifiers $u_i\colon\P\to\Fr$, $i\in I$, such that $\bigcap_{i\in I}\ker u_i$ is the identity relation on $P$. Thus, if there is a most general unifier for $\P$, its kernel is also the identity relation on $P$. But this would mean that $\P$ actually embeds into $\Fr$.
\end{proof}

\begin{proposition}\label{fact:: Medvedev and Levin}
Let $\V$ be one of the varieties $\V_{Med}$, $\V_{Lev_n}$ for $n\geq 2$. Then
\begin{enumerate}
\item $\V$ is SC,
\item there exists a $\V$-finitely presented unifiable algebra which does not embed into $\Fr$,
\item $\V$ does not have unitary (and hence projective)  unification.
\end{enumerate}
\end{proposition}

\begin{proof}
\mbox{}\\
\noindent (1) For $\V_{Med}$ it was proved by Prucnal \cite{Pruc76}. For $\V_{Lev_n}$ a small modification of Prucnal's proof works  \cite[Lemma 3.2]{Skv99}.

\noindent (2) It follows from Lemmas \ref{lem:: 2^2 finitely presented} and \ref{lem:: 2^2 not embeddable}.

\noindent (3) It follows from (1), (2) and Lemma \ref{lem:: no unitary unification}.
\end{proof}

In fact, point (3) of Proposition \ref{fact:: Medvedev and Levin} follows from \cite[Theorem 4.4]{Ghi99} and also from \cite[Lemmas 3,4]{Dzi06b}. It was shown there that every variety of Heyting algebras containing $\2^2\oplus\1$ cannot have unitary unification even without assuming SC. Indeed, then the algebra $\2^2$ does not have a most general unifier.

In the next section, based on Example \ref{exp:: Levin and Medvedev}, we will construct ASC$\setminus$SC varieties of closure algebras for which points (2) and (3) in Proposition \ref{fact:: Medvedev and Levin} will be also valid
\qed\end{example}

\section{More examples:  ASC for varieties of closure algebras and normal modal logics}\label{sec:: more examples}

In this section we will show that the varietal join $\V=\U \join \W$ of an SC variety of closure algebras $\U$ and a non-minimal variety of monadic algebras $\W$ is ASC$\setminus$SC.

By Corollary \ref{cor:: on finitely presented}, every $\V$-finitely presented unifiable algebra is isomorphic to a subalgebra of a power of $\Fr$. Therefore such finite $\P$ is a subalgebra of a power of some $\Fr(k)$. However in every known example of ASC$\setminus$SC variety $\V$ of modal algebras all \mbox{$\V$-finitely} presented unifiable algebras actually embed into $\Fr$. We claim that this fact is connected with the limitation of the techniques used so far, not with any intrinsic property of modal algebras. Indeed, a clue for it with varieties of Heyting algebras was already presented in Proposition \ref{fact:: Medvedev and Levin}. We will prove in Theorem \ref{thm:: modal companions of Medvedev and Levin} that if $\V$ is as in the previous paragraph and additionally $\U$ contains the algebra ${\sf B}(\2^2\oplus\1)$ shown in Figure \ref{fig::B(Lev_2^+)}, then there is a finite $\V$-finitely presented unifiable algebra which is not embeddable into $\Fr$. Moreover, $\V$ does not have unitary (and hence projective) unification.
For such $\U$ we can take minimal modal companions of Levin varieties from Example \ref{exp:: Levin and Medvedev}. To the best of our knowledge, these are the first found examples of ASC$\setminus$SC varieties of modal algebras without projective unification.

\subsection{Join of varieties of McKinsey algebras and of monadic algebras}

Let
\[
\mu(x)=\mBox\Diamond x\mto\Diamond\mBox x
\]
be the \emph{McKinsey} term. A modal algebra $\M$ is a \emph{McKinsey algebra} if
it satisfies the \emph{McKinsey identity}
\[
(\forall x)\;\mu(x)\=1.
\]
McKinsey algebras appeared in our investigations due to the fact that an ASC variety of closure algebras is SC iff it satisfies McKinsey identity. Moreover,
free algebras of finite rank for the varietal join of a variety of McKinsey algebras and of a variety of monadic algebras are products of McKinsey algebras and monadic algebras. These facts will be used to verify that a varietal join of an SC variety of closure algebras with a non-minimal variety of monadic algebras is ASC$\setminus$SC.
Now we will present their  proofs.

In what follows $\U$ will be a variety of McKinsey algebras, $\W$ will be a variety of monadic algebras, and $\V=\U\join\W={\sf V}(\U\cup \W)$ will be their varietal join. We will add subscripts in the notation of free algebras denoting varieties for which these algebras are free. For instance a free algebra for $\V$ of rank $\aleph_0$, a previously denoted $\Fr$, now will be denoted by $\Fr_\V$. Moreover, in this section $\2$ will again denote a two-element closure algebra.

Recall from Example \ref{exp:: monadic algebras} that we may put
\[
\Fr_\W(k)=\2^d\times\prod_{l=1}^{m}\R_l,
\]
where all $\R_l$ are not necessarily distinct finite simple monadic algebras with more than two elements, this means that they are in $\{\S_l\mid l\in \{2,3,\ldots\}\}$, and $d,m$ are some natural numbers. Let $w_1,\ldots,w_k$ be free generators of $\Fr_\W(k)$. Let us interpret them as mappings with the domain $\{0,\ldots,m\}$ and $w(0)\in 2^d$, $w(l)\in  R_l$ for $l\in\{1,\ldots,m\}$. Define
\[
\G_\W(k)=\prod_{l=1}^{m}\R_l.
\]
What we need to know about free generators in $\Fr_\W(k)$ is just the following fact.

\begin{lemma}\label{lem:: free monadic algebras}
For every index $l\in\{1,\ldots,m\}$ there exists $i\in\{1,\ldots,k\}$ such that \mbox{$w_i(l)\not\in\{0,1\}$}.
\end{lemma}

\begin{proof}
Indeed, otherwise the subalgebra of $\Fr_\W(k)$ generated by  $w_1,\ldots,w_k$ would be a subalgebra of $\2^d\times\prod_{j=1}^{l-1}\R_j\times\2\times\prod_{j=l+1}^m\R_{j}$. This would contradict the fact that free generators generate the whole algebra $\Fr_\W(k)$.
\end{proof}

\begin{proposition}\label{prop:: free algebras in McKinsey join monadic}
Let $\U$ be a nontrivial variety of McKinsey algebras, $\W$ be a variety of monadic algebras, and $\V=\U\join\W$  be its varietal join. Then
\[
\Fr_\V(k)\cong\Fr_{\U}(k)\times\G_\W(k).
\]
\end{proposition}


\begin{proof}
Let $u_1,\ldots,u_k$ be free generators of $\Fr_\U(k)$ and $w_1,\ldots,w_k$ be free generators of $\Fr_\W(k)$ interpreted as above.
For $i\in\{1,\ldots,k\}$ let $v_i=(u_i,w_i')$, where \mbox{$w'_i=w_i|_{\{1,\ldots, m\}}$}. We will prove that $v_1,\ldots,v_k$ are free generators in $\Fr_{\U}(k)\times\G_\W(k)$ for the variety $\V$. The verification of the universal mapping property may be split into the verification of the following two claims.

\begin{claim}
The elements $v_1,\ldots,v_k$ generate $\Fr_{\V}(k)\times\G_\W(k)$.
\end{claim}

\begin{proof}
The elements $u_1,\ldots, u_k$ generate $\Fr_\U(k)$ and $w'_1,\ldots, w'_k$ generate $\G_\W(k)$. Thus every element from $F_{\V}(k)\times G_\W(k)$ is of the form $(s(\bar u),t(\bar w'))$ for some terms $s(\bar x),t(\bar x)$. Our aim is to find a term $r(\bar x)$ such that $(s(\bar u),t(\bar w'))=r(\bar v)$.  Define a term
\[
m(\bar x)=\dsmbigmeet_{i=1}^k\mu(x_i).
\]
Since $\U$ satisfies McKinsey identity, in $\Fr_\U(k)$ we have $m(\bar u)=1$. Let us compute $m(\bar w')$ in $\G_\W(k)$. A routine verification shows that if $a$ is neither a top nor a bottom element in $\R_l$, then $\mu(a)=0$. Thus, by Lemma \ref{lem:: free monadic algebras}, for every $l\in\{1,\ldots,m\}$ there is $i$ such that in $\R_l$ we have
$\mu(w_i(l))=0$.
Hence in $\G_\W(k)$ we have $m(\bar w')=0$. Now for $r(\bar x)=(m(\bar x)\mmeet s(\bar x))\mjoin(\mneg m(\bar x)\mmeet t(\bar x))$ we compute
\begin{multline*}
r(\bar v)=((m(\bar u)\mmeet s(\bar u))\mjoin(\mneg m(\bar u)\mmeet t(\bar u),(m(\bar w')\mmeet s(\bar w'))\mjoin(\mneg m(\bar w')\mmeet t(\bar w'))\\
=((1\mmeet s(\bar u))\mjoin(0\mmeet t(\bar u),(0\mmeet s(\bar w'))\mjoin(1\mmeet t(\bar w'))=(s(\bar u),t(\bar w')).
\end{multline*}
\end{proof}

\begin{claim}
For every algebra $\M\in\V$ and every mapping $f\colon\{v_1,\ldots,v_k\}\to M$ there is a homomorphism $\bar f\colon\Fr_{\U}(k)\times\G_\W(k)\to \M$ extending $f$.
\end{claim}

\begin{proof}
First observe that we do not have to verify the assertion for all $\M\in\V$. It is enough to show this for generators of $\V$, see e.g. \cite[Proposition 4.8.9]{Kra99}. We will do it for algebras from $\U\cup\W$.

The case when $\M\in \U$ is easy. Then for $\bar f$ we just take the composition of the first projection of $\Fr_{\U}(k)\times\G_\W(k)$ with the homomorphism from $\Fr_\U(k)$ into $\M$ extending the mapping given by $u_i\mapsto f(v_i)$ for $i\in\{1,\ldots,k\}$.

Let us move to the case when $\M\in \W$.  Since we assumed that $\U$ is nontrivial, $\2^d\in\U$. (Actually $\2^d$ is  free for $\sf V(\2)$ of rank $k$ and $d=2^k$.) Thus there is a homomorphism $g\colon\Fr_\U(k)\to \2^{d}$  such that $g(u_i)=w_i(0)$ for $i\in\{1,\ldots,k\}$ (recall that $w_i$ are free generators in $\Fr_\W(k)$). Let
\[
g'\colon\Fr_{\U}(k)\times\G_\W(k)\to\Fr_\W(k);\; (s,t)\mapsto (g(s),t).
\]
and $h\colon \Fr_\W(k)\to\M$ be a homomorphism such that $h(w_i)=f(v_i)$ for $i\in\{1,\ldots,k\}$. Then $h\circ g'\colon \Fr_{\U}(k)\times\G_\W(k)\to\M$ is a homomorphism extending $f$.
\end{proof}
\end{proof}

In \cite{DW12} it was proved that a variety of S4.3-algebras is SC iff it satisfies McKinsey identity. However the proof presented there actually uses only the fact that a variety of closure algebras under consideration is ASC. Let us formulate the reasoning in algebraic terms. The following lemma will be used several times.

\begin{lemma}\label{lem:: homo McKinsey to 2}
For every closure McKinsey algebra $\M$ and every open element $a\in M$ which does not equal 0 there is a homomorphism $h\colon \M\to \2$ such that $h(a)=1$.
\end{lemma}

\begin{proof}
Let $\theta$ be a congruence of $\M$ generated by $(a,1)$.
Since $a$ is open, $0\not\in 1/\theta$, and therefore $\theta<M^2$.
By Zorn lemma, $\theta$ can be extended to a maximal congruence $\theta'$ properly contained in $M^2$. Then $\N  = \M/\theta'$ is simple. For closure algebras being simple is equivalent to having exactly two open elements 0 and 1. The computation shows that in $\N$ we have $\mu(c)=0$ for every element $c \in N-\{0,1\}$. Thus the inequality $N\neq\{0,1\}$ would contradict the satisfaction of McKinsey identity by $\N$ (see also Proposition \ref{prop:: splitting McKinsey}). Therefore $\N \cong \2$.
\end{proof}

Recall that the four-element simple closure algebras, depicted in Figure \ref{fig:: 2 and S_2}, was denoted by $\S_2$. Note that $\sf V(\S_2)$ and the variety of McKinsey algebras is a splitting pair for the lattice of varieties of closure algebras.


\begin{proposition}[\protect{\cite[Example III.3.9]{Blo76}, \cite[Example IV.5.4]{Rau79}}]\label{prop:: splitting McKinsey}
Let $\U$ be a variety of closure algebras. Then
$\S_2\not\in\U$  if and only if  $\U$ satisfies McKinsey identity.
\end{proposition}


\begin{lemma}\label{lem:: SC -> McKinsey}
Every SC variety of closure algebras satisfies McKinsey identity.
\end{lemma}

\begin{proof}
On the contrary, assume that $\U$ does not satisfy McKinsey identity. Then, by Proposition \ref{prop:: splitting McKinsey}, $\S_2\in\U$.
Let
\[
q=(\forall x)[\Diamond x   \mmeet  \Diamond \mneg x \=1 \to 0\=1].
\]
We have $\S_2\not\models q$, and hence $\U\not\models q$.
But $q$ is $\sf V(\2)$-passive. Since $\Fr_{{\sf V}(2)}$ is a homomorphic image of $\Fr_\U$, the quasi-identity $q$ is also $\U$-passive, and therefore it holds in $\Fr_\U$. Thus $\U$ is not SC.
\end{proof}

\begin{proposition}\label{prop:: SC <-> McK}
Let $\U$ be an ASC variety of closure algebras. Then the following conditions are equivalent:
\begin{enumerate}
\item $\U$ is SC,
\item $\U$ satisfies McKinsey identity,
\item $\S_2\not\in\U$.
\end{enumerate}
\end{proposition}

\begin{proof}
The equivalence (2)$\Leftrightarrow$(3) follows from Proposition \ref{prop:: splitting McKinsey}, the implication  (1)$\Rightarrow$(2) is given by Lemma \ref{lem:: SC -> McKinsey} and the implication (2)$\Rightarrow$(1)  follows from Corollary \ref{cor:: ASC <-> SC}, and Lemma \ref{lem:: homo McKinsey to 2}.
\end{proof}

Now we may proceed to the main result of this section.

\begin{theorem}\label{thm:: ASC for McK join Monadic}
Let $\U$ be an SC variety of closure algebras and $\W$ be a non-minimal variety of monadic algebras.
Then the varietal join $\U=\V \join\W$ is ASC$\,\setminus$SC.
\end{theorem}

\begin{proof}
In the case when $\U$ is trivial we have $\V=\U$ and the statement of the theorem was verified in Example \ref{exp:: monadic algebras}. So we assume that $\U$ is nontrivial.
Let us start by proving the following fact.

\begin{claim}
The algebra $\Fr_\U(k)$ embeds into $\Fr_\V(k)$ for every natural number $k$. In particular, ${\sf Q}(\Fr_\U)\leq{\sf Q}(\Fr_\V)$.
\end{claim}

\begin{proof}
By Lemma \ref{lem:: homo McKinsey to 2}, there is a homomorphism $h\colon\Fr_\U(k)\to\G_\W(k)$ with the image isomorphic to $\2$. Then the homomorphism
\[
g\colon\Fr_\U(k)\to\Fr_\U(k)\times\mathbf G_\W(k);\; t\mapsto(t,h(t))
\]
embeds $\Fr_\U(k)$ into $\Fr_\U(k)\times\G_\W(k)$. By Propositions \ref{prop:: free algebras in McKinsey join monadic} and \ref{prop:: SC <-> McK}, the later algebra is isomorphic to $\Fr_\V(k)$.
\end{proof}

In order to verify ASC for $\V$ let us check the condition (3') from Corollary \ref{cor:: about C}.
Since $\V$ is congruence distributive, every subdirectly irreducible algebra $\S$ from $\V$ is in $\U\cup\W$ \cite[Corollary 4.2]{Jon67}.

\noindent{\it Case when $\S\in\U$.} By the assumption that $\U$ is SC, $\S\in{\sf Q(\Fr_\U)}$. Hence $\S\times\2\in{\sf Q(\Fr_\U)}$, and by Claim, $\S\times\2\in{\sf Q(\Fr_\V)}$.

\noindent{\it Case when $\S\in\W$.} As explained in Example \ref{exp:: monadic algebras}, for large enough $k$, $\S$ is a proper direct factor of $\Fr_\W(k)$. Thus, by  Proposition \ref{prop:: free algebras in McKinsey join monadic}, $\S$ is a proper direct factor of $\Fr_\V(k)$, and  hence $\S\times\2$ embeds into $\Fr_\V(k)$.

Finally note that, since $\W$ is non-minimal, $\S_2\in\V$. Thus by Proposition \ref{prop:: SC <-> McK}, $\V$ is not SC.
\end{proof}

\subsection{ASC without projective unification; modal companions of Levin and Medvedev varieties}.

Let us briefly review the translation of intuitionistic logics into transitive reflexive normal modal logics in algebraic terms.  Open elements of every closure algebra $\M$ form the Heyting algebra $\sf O(\M)$ with the order inherited from $\M$. Moreover, if $\W$ is a variety of closure algebras, then $\sf O(\W)$ is a variety of Heyting algebras ($\sf O$ is treated here as a class operator). The following fact was proved in \cite[Section 1]{MT46}, see also \cite[Chapter 1]{Blo76} and \cite[Theorem 2.2]{BD75}.

\begin{proposition}\label{prop:: McKinsey-Tarski thm}
For every Heyting algebra $\H$ there is a closure algebra ${\sf B}(\H)$ such that
\begin{enumerate}
\item ${\sf OB}(\H)=\H$;
\item for every closure algebra $\M$, if
$\H\leq{\sf O}(\M)$, then ${\sf B}(\H)$ is isomorphic to the subalgebra of $\M$ generated by $H$,
\end{enumerate}
\end{proposition}

The algebra ${\sf B}(\H)$ is called the \emph{free Boolean extension of} $\H$.

For each variety $\mathcal Y$ of Heyting algebras there is a variety of closure algebras $\mathcal Y'$ such that ${\sf O}(\mathcal Y')=\mathcal Y$. Every such $\mathcal Y'$ is called a \emph{modal companion} of $\mathcal Y$. In general, there are many modal companions for a given variety of Heyting algebras. For instance the variety of monadic algebras is the greatest modal companion, and the variety $\sf V(\2)$ is the smallest modal companion of the variety of Boolean algebras treated as Heyting algebras. The smallest modal companion of a variety $\mathcal Y$ of Heyting algebras is given by
\[
\sigma(\mathcal Y)={\sf V}\{{\sf B}(\H)\mid \H\in\mathcal Y\}
\]
The importance of $\sigma$ and $\sf O$ class operators follows from Blok-Esakia theorem. It states that they  are mutually inverse lattice isomorphisms between the subvariety lattice of the variety of Heyting algebras and the subvariety lattice of the variety of Grzegorczyk algebras \cite[Theorem III.7.10]{Blo76}\cite[Theorem 9.66]{CZ97} \cite{ZW12} \cite{Esa79b}.

As in Example \ref{exp:: Levin and Medvedev} we would like to find a $\V$-finitely presented algebra that does not embed into $\Fr_\V$. But this time we cannot take $\2^2$. Indeed, Proposition \ref{prop:: free algebras in McKinsey join monadic} shows that this algebra in fact embeds into $\Fr_\V$. But a bit more complicated algebra will not embed. Let $\4$ be a four-element closure algebra depicted in Figure \ref{fig::4}
\begin{figure}[h]
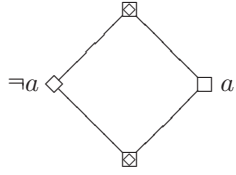

\[
\xy
(10,20)*{\DiamondBox}="AC";
(0,9.8)*{\Diamond}="A";  (20,10)*{\mBox}="C";
(10,0)*{\DiamondBox}="ZERO";
(-4.1,10.1)*{\mneg a};
(23,10)*{a};
(9.07,1.07); (.6,9.54) **\dir{-}; 
(10.93,1.07);  (19.1,9.24) **\dir{-}; 
(9.07,19.22); (.60,10.75) **\dir{-}; 
(19.07,11.07); (10.9,19.24) **\dir{-};  
\endxy
\]
\caption{The closure algebra ${\4}$.}
\label{fig::4}
\end{figure}
Note that $\sf V(\4)$ and the variety of monadic algebras is a splitting pair for the lattice of varieties of closure algebras \cite[Theorem 5.5]{BD75}.  We will not need this fact in full strength, but only the observation that $\4$ embeds into ${\sf B}(\2^2\oplus\1)$.

\begin{lemma}\label{lem:: 4^2 is finitely presented}
Let $\V$ be a variety of closure algebras containing $\4$. Then $\4^2$ is a $\V$-finitely presented unifiable algebra.
\end{lemma}

\begin{proof}
Fact \ref{fact::finite->fp} yields that $\4^2$ is $\V$-finitely presented. With some work one may also show that $\4^2\cong\P_{\varphi(x)}$, where 
\[
\varphi= (\mBox\Diamond\mBox v\=\Diamond\mBox v) \meet (\Diamond\mBox v \mmeet v\= \mBox v) \meet (\Diamond\mBox v \mjoin v\= \Diamond v)
 \]
Moreover, there is a homomorphism from $\4$ onto $\2$. Hence $\4^2$ is unifiable.
\end{proof}

In our considerations the algebra ${\sf B}(\2^2\oplus\1)$ depicted in Figure \ref{fig::B(Lev_2^+)}, which is the free Boolean extension of the Heyting algebra $\2^2\oplus\1$, plays a crucial role. Note that ${\sf B}(\2^2\oplus\1)$ is a modal algebra which is dual to the ordered set $\Lev_2$ treated as a modal frame.
\begin{figure}[h]
\[
\xy
(10,30)*{\DiamondBox}="ONE";
(0.01,19.81)*{\Diamond}="AB"; (10,20)*{\mBox}="AC"; (20.01,19.81)*{\Diamond}="BC";
(0,10)*{\mBox}="A"; (10.01,9.78)*{\Diamond}="B"; (20,10)*{\mBox}="C";
(10,0)*{\DiamondBox}="ZERO";
(-3,10)*{a}; (13,10)*{c}; (23,10.43)*{b};
(9.07,1.07); (.9,9.24) **\dir{-}; 
(10,1.1); (10,9) **\dir{-}; 
(10.93,1.07);  (19.1,9.24) **\dir{-}; 
(0,11.1); (0,19.0) **\dir{-}; 
(9.07,19.22); (.90,11.05) **\dir{-}; 
(9.40,10.69); (0.60,19.59) **\dir{-}; 
(19.39,19.51); (10.60,10.72) **\dir{-}; 
(19.07,11.07); (10.9,19.24) **\dir{-};  
(20,11.1); (20,19.0) **\dir{-}; 
(9.07,29.22); (.60,20.75) **\dir{-}; 
(10,21.1); (10,29.2) **\dir{-}; 
(19.42,20.73); (10.95,29.20) **\dir{-}; 
\endxy
\]
\caption{The closure algebra $\sf B(\2^2\oplus\1$).}
\label{fig::B(Lev_2^+)}
\end{figure}

\begin{lemma}\label{lem:: 4^2 is not embeddable}
Let $\U$ be a variety of closure McKinsey algebras containing ${\sf B}(\2^2\oplus\1$), $\W$ be a variety of monadic algebras and $\V=\U\join\W$ be their varietal join. Then $\4^2$ does not embed into $\Fr_\V$.
\end{lemma}

\begin{proof}
In order to obtain a contradiction, assume that $\4^2$ embeds into $\Fr$. Then, since $\4^2$ is finite, it embeds into $\Fr_\V(k)$ for some natural number $k$. Recall that, by Theorem \ref{thm:: ASC for McK join Monadic}, $\Fr_\V(k)\cong\Fr_\U(k)\times\G_\W(k)$.
The congruence lattice of $\4^2$ is isomorphic to a product of two three-element chains. Let $\rho$ be the congruence from the middle of this square, i.e., $\rho$ is generated by $((a,a),(1,1))$. Now, if $\R$ is a simple algebra in $\W$ and $h\colon \4^2\to\R$ is a homomorphism, then $\rho\leq\ker(h)$. Thus for every homomorphism $h\colon\4^2\to\G_\W(k)$ we also have $\rho\leq\ker(h)$. Moreover, $\rho\cap \alpha=\id_{4^2}$ iff $\alpha=\id_{4^2}$ for every congruence $\alpha$ of $\4^2$. These facts yield that $\4^2$ must embed into $\Fr_\U(k)$. Since $\2^2\leq\4^2$, $\2^2$ also embeds into $\Fr_\U(k)$.

The rest of the proof is very similar to the proof of Lemma \ref{lem:: 2^2 not embeddable}.
The generator $t$ of an isomorphic image of $\2^2$ in $\Fr_\V(k)$ satisfies
\[
0<t,\mneg t<1\quad\text{and}\quad \mBox t =t,\quad \mBox \mneg t =\mneg t.
\]
Therefore, by Lemma \ref{lem:: homo McKinsey to 2}, there exists a homomorphism $f\colon\Fr_\U\to \2^2$ such that $f(t)=(1,0)$ and $f(\mneg t)=(0,1)$. Let $a=\{(1,0)\}$, $b=\{(0,1)\}$ be the open atoms and $a\mjoin b =\{(0,1),(1,0)\}$ be the open coatom in ${\sf B}(\2^2\oplus\1)$. Let $g\colon V\to B(2^2\oplus 1)$ be a mapping given by
\[
g(v)=
\begin{cases}
a &\text{ if } f(v)=(1,0)\\
b &\text{ if } f(v)=(0,1)\\
a\mjoin b&\text{ if } f(v)=(1,1)\\
0 &\text{ if } f(v)=(0,0)
\end{cases},
\]
and  $\bar g\colon\Fr\to{\sf B}(\2^2\oplus\1)$ be the homomorphic extension of $g$. Let $h\colon{\sf B}(\2^2\oplus\1)\to\2^2$ be the surjective homomorphism such that $h(a)=(1,0)$ and $h(b)=(0,1)$. Then $f|_V=h\circ g$. Thus, by the universal mapping property of $\Fr_\V(k)$, $f=h\circ\bar g$. We get that $\bar g(t)\in h^{-1}((1,0))=\{a,a\mjoin c\}$ and $\bar g(\mneg t)\in h^{-1}((0,1))=\{b,b\mjoin c\}$, where $c$ is the third, closed atom of ${\sf B}(\2^2\oplus\1)$. Therefore, since $\bar g$ maps open elements onto open element,  $\bar g( t)= a$ and $\bar g( \mneg t)= b$. Now we compute in ${\sf B}(\2^2\oplus\1)$
\[
1=\bar g(1)=\bar g(t\mjoin\mneg t)=\bar g(t)\mjoin\mneg \bar g(t)=a\mjoin b<1,
\]
and reach a contradiction.

We would like to finish this proof with one technical remark. In the earlier version of this paper we dealt with $\sigma(\mathcal Y)$, where $\mathcal Y$ is a variety of Heyting algebras containing $\2^2\oplus\1$, instead of $\U$. One would then wish to use Lemma \ref{lem:: 2^2 not embeddable} instead of repeating the whole argumentation as we did here. But it does not give a correct proof. It follows from the fact that in general $\Fr_{\mathcal Y}(k)$ is not isomorphic to ${\sf O}(\Fr_{\sigma(\mathcal Y)}(k)$). Actually, $\Fr_{\mathcal Y}(k)$ only embeds into ${\sf O}(\Fr_{\sigma(\mathcal Y)}(k)$) and this embedding is proper even in the simple case when $\mathcal Y=\V_{Lev_2}$.
\end{proof}

\begin{theorem}\label{thm:: modal companions of Medvedev and Levin}
Let $\U$ be an SC variety of closure algebras containing ${\sf B}(\2^2\oplus\1)$, for instance any of $\sigma(\V_{Lev_2}),\sigma(\V_{Lev_3}),\ldots,
\sigma(V_{Med})$, $\W$ be a non-minimal variety of monadic algebras and $\V=\U\join\W$ be their varietal join. Then
\begin{enumerate}
\item $\V$ is ASC$\,\setminus$SC,
\item there exists a $\V$-finitely presented unifiable algebra which does not embed into $\Fr$,
\item $\V$ does not have unitary (and hence projective)  unification.
\end{enumerate}
\end{theorem}

\begin{proof}
\mbox{}\\
\noindent (1) It is a consequence of  Theorem  \ref{thm:: ASC for McK join Monadic}. Notice that all varieties among $\sigma(\V_{Lev_2})$, $\sigma(\V_{Lev_3})$, $\ldots$,
$\sigma(V_{Med})$ are SC. Indeed, it follows from Proposition \ref{fact:: Medvedev and Levin} and the fact that $\sigma$ operator preserves SC \cite[Theorem 5.4.7]{Ryb97}.

\noindent (2) Since, $\4$ embeds into ${\sf B}(\2^2\oplus\1)$,  $\4^2\in\V$. Thus it follows from Lemmas \ref{lem:: 4^2 is finitely presented} and \ref{lem:: 4^2 is not embeddable}.

\noindent (3) It follows from (1), (2) and Lemma \ref{lem:: no unitary unification}. This point was also proved in \cite[Lemmas 3,4]{Dzi06b}, where it was shown that every variety of closure algebras, not necessarily SC,  containing ${\sf B}(\2^2\oplus\1)$ cannot have unitary unification. Again, then $\2^2$ is a $\V$-finitely presented algebra without a most general unifier.
\end{proof}

Let us note that the mapping $(\U,\V)\mapsto\U\join\V$, where $\U$ and $\V$ are nontrivial varieties satisfying the condition from Theorem \ref{thm:: ASC for McK join Monadic}, is injective. To see this assume that $\V=\U^0\join\W^0=\U^1\join\W^1$, where $\U^0,\U^1$ and $\W^0,\W^1$ satisfies the same conditions as $\U$ and $\W$ respectively. By \cite[Corollary 4.2]{Jon67}, $$\V_{SI}=\U^0_{SI}\cup\W^0_{SI}=\U^1_{SI}\cup\W^1_{SI}.$$ Since $\U^i$ are varieties of McKinsey algebras and $\W^j$ are varieties of Monadic algebras, for $i,j\in \{0,1\}$, we have
\[
\U^i\cap\V^j=\{\2\}.
\]
This yields the equations $\U^0_{SI}=\U^1_{SI}$ and $\V^0_{SI}=\V^1_{SI}$. Hence $\U^0=\U^1$ and $\V^0=\V^1$.

Therefore, there are at least as many ASC$\setminus$SC varieties of closure algebras as there are SC varieties of closure algebras. However, we do not know exactly how many of them there are. We know that there
are at least $\aleph_0$ (the number of Levin varieties) and at most $\mathfrak{c}$ (the number of varieties of closure algebras) members of both families. Note that there are continuum many ASC varieties of modal algebras \cite[Corollary 14]{Dzi06}.

\begin{problem}
How many SC and ASC varieties of closure algebras are there?
\end{problem}

\subsection*{Acknowledgement}
We wish to thank Alex Citkin who suggested that ASC$\setminus$SC property might be detected among ASC varieties by possessing some simple algebra. We confirmed this for varieties of closure algebras in Proposition \ref{prop:: SC <-> McK}.

\bibliographystyle{plain}
\bibliography{asc11082014}

\begin{thebibliography}{10}

\bibitem{AJN91}
Hajnal Andr\'eka, Bjarni J\'onsson, and Istv\'an N\'emeti.
\newblock Free algebras in discriminator varieties.
\newblock {\em Algebra Universalis}, 28(3):401--447, 1991.

\bibitem{Bas58}
Hyman Bass.
\newblock Finite monadic algebras.
\newblock {\em Proc. Amer. Math. Soc.}, 9:258--268, 1958.

\bibitem{Ber88}
Clifford Bergman.
\newblock Structural completeness in algebra and logic.
\newblock In {\em Algebraic logic ({B}udapest, 1988)}, volume~54 of {\em
  Colloq. Math. Soc. J\'anos Bolyai}, pages 59--73. North-Holland, Amsterdam,
  1991.

\bibitem{Blo76}
Willem~J. Blok.
\newblock {\em Varieties of interior algebras}.
\newblock PhD thesis, University of Amsterdam, 1976.
\newblock
  URL=\url{http://www.illc.uva.nl/Research/Dissertations/HDS-01-Wim_Blok.text.pdf}.

\bibitem{BD75}
Willem~J. Blok and Philip Dwinger.
\newblock Equational classes of closure algebras. {I}.
\newblock {\em Nederl. Akad. Wetensch. Proc. Ser. A {\bf 78}=Indag. Math.},
  37:189--198, 1975.

\bibitem{BP82}
Willem~J. Blok and Don Pigozzi.
\newblock On the structure of varieties with equationally definable principal
  congruences. {I}.
\newblock {\em Algebra Universalis}, 15(2):195--227, 1982.

\bibitem{BP89}
Willem~J. Blok and Don Pigozzi.
\newblock Algebraizable logics.
\newblock {\em Mem. Amer. Math. Soc.}, 77(396):vi+78, 1989.

\bibitem{BP88}
Willem~J. Blok and Don Pigozzi.
\newblock Local deduction theorems in algebraic logic.
\newblock In {\em Algebraic logic ({B}udapest, 1988)}, volume~54 of {\em
  Colloq. Math. Soc. J\'anos Bolyai}, pages 75--109. North-Holland, Amsterdam,
  1991.

\bibitem{BP01}
Willem~J. Blok and Don Pigozzi.
\newblock Abstract algebraic logic and the deduction theorem, 2001.
\newblock Manuscript available at
  \url{http://orion.math.iastate.edu/dpigozzi/}.

\bibitem{BvA02}
Willem~J. Blok and Clint~J. van Alten.
\newblock The finite embeddability property for residuated lattices, pocrims
  and bck-algebras.
\newblock {\em Algebra Universalis}, 48(3):253--271, 2002.

\bibitem{Bul66}
Robert~A. Bull.
\newblock That all normal extensions of {${\rm S}4.3$} have the finite model
  property.
\newblock {\em Z. Math. Logik Grundlagen Math.}, 12:341--344, 1966.

\bibitem{Bur92}
Stanley Burris.
\newblock Discriminator varieties and symbolic computation.
\newblock {\em J. Symbolic Comput.}, 13(2):175--207, 1992.

\bibitem{BS81}
Stanley Burris and H.~P. Sankappanavar.
\newblock {\em A course in universal algebra}, volume~78 of {\em Graduate Texts
  in Mathematics}.
\newblock Springer-Verlag, New York, 1981.
\newblock The Millennium Edition is available at
  \url{http://www.math.uwaterloo.ca/~snburris/htdocs/ualg.html}.

\bibitem{CSV14}
Miguel~A. Campercholi, Micha\l~M. Stronkowski, and Diego~J. Vaggione.
\newblock On structural completeness vs almost structural completeness problem:
  {A} discriminator varieties case study.
\newblock arXiv:1407.0175, 2014.

\bibitem{CZ97}
Alexander Chagrov and Michael Zakharyaschev.
\newblock {\em Modal logic}, volume~35 of {\em Oxford Logic Guides}.
\newblock The Clarendon Press Oxford University Press, New York, 1997.
\newblock Oxford Science Publications.

\bibitem{CDM00}
Roberto L.~O. Cignoli, Itala M.~L. D'Ottaviano, and Daniele Mundici.
\newblock {\em Algebraic foundations of many-valued reasoning}, volume~7 of
  {\em Trends in Logic}.
\newblock Kluwer Academic Publishers, Dordrecht, 2000.

\bibitem{CM09}
Petr Cintula and George Metcalfe.
\newblock Structural completeness in fuzzy logics.
\newblock {\em Notre Dame J. Form. Log.}, 50(2):153--182, 2009.

\bibitem{Cze01}
Janusz Czelakowski.
\newblock {\em Protoalgebraic logics}, volume~10 of {\em Trends in
  Logic---Studia Logica Library}.
\newblock Kluwer Academic Publishers, Dordrecht, 2001.

\bibitem{CzD92}
Janusz Czelakowski and Wies{\l}aw Dziobiak.
\newblock A single quasi-identity for a quasivariety with the fraser-horn
  property.
\newblock {\em Algebra Universalis}, 29(1):10--15, 1992.

\bibitem{DP02}
B.~A. Davey and H.~A. Priestley.
\newblock {\em Introduction to lattices and order}.
\newblock Cambridge University Press, New York, second edition, 2002.

\bibitem{DH01}
J.~Michael Dunn and Gary~M. Hardegree.
\newblock {\em Algebraic methods in philosophical logic}, volume~41 of {\em
  Oxford Logic Guides}.
\newblock The Clarendon Press Oxford University Press, New York, 2001.
\newblock Oxford Science Publications.

\bibitem{Dzi06b}
Wojciech Dzik.
\newblock Splittings of lattices of theories and unification types.
\newblock In {\em Contributions to general algebra. 17}, pages 71--81. Heyn,
  Klagenfurt, 2006.

\bibitem{Dzi06}
Wojciech Dzik.
\newblock Transparent unifiers in modal logics with self-conjugate operators.
\newblock {\em Bull. Sect. Logic Univ. \L\'od\'z}, 35(2-3):73--83, 2006.

\bibitem{Dzi08}
Wojciech Dzik.
\newblock Unification in some substructural logics of {BL}-algebras and hoops.
\newblock {\em Rep. Math. Logic}, 43:73--83, 2008.

\bibitem{Dzi11}
Wojciech Dzik.
\newblock Remarks on projective unifiers.
\newblock {\em Bull. Sect. Logic Univ. \L\'od\'z}, 40(1-2):37--46, 2011.

\bibitem{DW14}
Wojciech Dzik and Piotr Wojtylak.
\newblock Modal consequence relations extending {S}4.3. an application of
  projective unification.
\newblock {\em Notre Dame J. Form. Log.}
\newblock To appear.

\bibitem{DW12}
Wojciech Dzik and Piotr Wojtylak.
\newblock Projective unification in modal logic.
\newblock {\em Log. J. IGPL}, 20(1):121--153, 2012.

\bibitem{Esa79b}
Leo~L. Esakia.
\newblock On the theory of modal and superintuitionistic systems.
\newblock In {\em Logical inference ({M}oscow, 1974)}, pages 147--172.
  ``Nauka'', Moscow, 1979.

\bibitem{FR09}
Josep~M. Font and Ramon Jansana.
\newblock {\em A general algebraic semantics for sentential logics}, volume~7
  of {\em Lecture Notes in Logic}.
\newblock Association for Symbolic Logic, Berlin, second edition, 2009.
\newblock URL=\url{http://projecteuclid.org/euclid.lnl/1235416965}.

\bibitem{FJP03}
Josep~M. Font, Ramon Jansana, and Don Pigozzi.
\newblock A survey of abstract algebraic logic.
\newblock {\em Studia Logica}, 74(1-2):13--97, 2003.
\newblock Abstract algebraic logic, Part II (Barcelona, 1997).

\bibitem{Fre79}
Ralph Freese.
\newblock The variety of modular lattices is not generated by its finite
  members.
\newblock {\em Trans. Amer. Math. Soc.}, 255:277--300, 1979.

\bibitem{Ghi97}
Silvio Ghilardi.
\newblock Unification through projectivity.
\newblock {\em J. Logic Comput.}, 7(6):733--752, 1997.

\bibitem{Ghi99}
Silvio Ghilardi.
\newblock Unification in intuitionistic logic.
\newblock {\em J. Symbolic Logic}, 64(2):859--880, 1999.

\bibitem{Ghi00}
Silvio Ghilardi.
\newblock Best solving modal equations.
\newblock {\em Ann. Pure Appl. Logic}, 102(3):183--198, 2000.

\bibitem{Gor98}
Viktor~A. Gorbunov.
\newblock {\em Algebraic theory of quasivarieties}.
\newblock Siberian School of Algebra and Logic. Consultants Bureau, New York,
  1998.
\newblock Translated from the Russian.

\bibitem{Hal56}
Paul~R. Halmos.
\newblock Algebraic logic. {I}. {M}onadic {B}oolean algebras.
\newblock {\em Compositio Math.}, 12:217--249, 1956.

\bibitem{Hal59}
Paul~R. Halmos.
\newblock Free monadic algebras.
\newblock {\em Proc. Amer. Math. Soc.}, 10:219--227, 1959.

\bibitem{HMT85}
Leon Henkin, J.~Donald Monk, and Alfred Tarski.
\newblock {\em Cylindric algebras. {P}art {I} and {II}}, volume~64 of {\em
  Studies in Logic and the Foundations of Mathematics}.
\newblock North-Holland Publishing Co., Amsterdam, 1985.

\bibitem{HB85}
Leslie Hogben and Clifford Bergman.
\newblock Deductive varieties of modules and universal algebras.
\newblock {\em Trans. Amer. Math. Soc.}, 289(1):303--320, 1985.

\bibitem{Iem14b}
Rosalie Iemhoff.
\newblock On rules.
\newblock {\em J. Philos. Logic}.
\newblock To appear.

\bibitem{Iem14}
Rosalie Iemhoff.
\newblock Unification in transitive reflexive modal logics.
\newblock {\em Notre Dame J. Form. Log.}
\newblock To appear.

\bibitem{Igo74}
V.~I. Igo{\v{s}}in.
\newblock Quasivarieties of lattices.
\newblock {\em Mat. Zametki}, 16:49--56, 1974.

\bibitem{Jer05}
Emil Je{\v{r}}{\'a}bek.
\newblock Admissible rules of modal logics.
\newblock {\em J. Logic Comput.}, 15(4):411--431, 2005.

\bibitem{Jon67}
Bjarni J{\'o}nsson.
\newblock Algebras whose congruence lattices are distributive.
\newblock {\em Math. Scand.}, 21:110--121, 1967.

\bibitem{KQ76}
Joel Kagan and Robert~W. Quackenbush.
\newblock Monadic algebras.
\newblock {\em Rep. Math. Logic}, 7:53--61, 1976.

\bibitem{Kow14}
Tomasz Kowalski.
\newblock {BCK} is not structurally complete.
\newblock {\em Notre Dame J. Form. Log.}, 55(2):197--204, 2014.

\bibitem{Kra99}
Marcus Kracht.
\newblock {\em Tools and techniques in modal logic}, volume 142 of {\em Studies
  in Logic and the Foundations of Mathematics}.
\newblock North-Holland Publishing Co., Amsterdam, 1999.

\bibitem{Lev69}
Leonid~A. Levin.
\newblock Some syntactic theorems on {J}u. {T}. {M}edvedev's calculus of finite
  problems.
\newblock {\em Dokl. Akad. Nauk SSSR}, 185:32--33, 1969.

\bibitem{Mad06}
Roger~D. Maddux.
\newblock {\em Relation algebras}, volume 150 of {\em Studies in Logic and the
  Foundations of Mathematics}.
\newblock Elsevier B. V., Amsterdam, 2006.

\bibitem{Mal70}
Anatoly~I. Mal'cev.
\newblock {\em Algebraic systems}.
\newblock Springer-Verlag, New York-Heidelberg, 1973.

\bibitem{McK43}
J.~C.~C. McKinsey.
\newblock The decision problem for some classes of sentences without
  quantifiers.
\newblock {\em J. Symbolic Logic}, 8:61--76, 1943.

\bibitem{MT46}
J.~C.~C. McKinsey and Alfred Tarski.
\newblock On closed elements in closure algebras.
\newblock {\em Ann. of Math. (2)}, 47:122--162, 1946.

\bibitem{MT48}
J.~C.~C. McKinsey and Alfred Tarski.
\newblock Some theorems about the sentential calculi of {L}ewis and {H}eyting.
\newblock {\em J. Symbolic Logic}, 13:1--15, 1948.

\bibitem{Med62}
Ju.~T. Medvedev.
\newblock Finitive problems.
\newblock {\em Dokl. Akad. Nauk SSSR}, 142:1015--1018, 1962.

\bibitem{med66}
Ju.~T. Medvedev.
\newblock Interpretation of logical formulas by means of finite problems.
\newblock {\em Soviet Math. Dokl.}, 7:857--860, 1966.

\bibitem{MR13}
George Metcalfe and Christoph R{\"o}thlisberger.
\newblock Admissibility in finitely generated quasivarieties.
\newblock {\em Log. Methods Comput. Sci.}, 9(2):2:09, 19, 2013.

\bibitem{ORV08}
J.~S. Olson, J.~G. Raftery, and C.~J. van Alten.
\newblock Structural completeness in substructural logics.
\newblock {\em Log. J. IGPL}, 16(5):455--495, 2008.

\bibitem{OR07}
Jeffrey~S. Olson and James~G. Raftery.
\newblock Positive {S}ugihara monoids.
\newblock {\em Algebra Universalis}, 57(1):75--99, 2007.

\bibitem{Pog71}
Witold~A. Pogorzelski.
\newblock Structural completeness of the propositional calculus.
\newblock {\em Bull. Acad. Polon. Sci. S\'er. Sci. Math. Astronom. Phys.},
  19:349--351, 1971.

\bibitem{PW08}
Witold~A. Pogorzelski and Piotr Wojtylak.
\newblock {\em Completeness theory for propositional logics}.
\newblock Studies in Universal Logic. Birkh\"auser Verlag, Basel, 2008.

\bibitem{Pruc76}
Tadeusz Prucnal.
\newblock Structural completeness of {M}edvedev's propositional calculus.
\newblock {\em Rep. Math. Logic}, 6:103--105, 1976.

\bibitem{Raf14}
James~G. Raftery.
\newblock Admissible rules and the {L}eibniz hierarchy.
\newblock {\em Notre Dame J. Form. Log.}
\newblock To appear.

\bibitem{Raf11}
James~G. Raftery.
\newblock Contextual deduction theorems.
\newblock {\em Studia Logica}, 99(1-3):279--319, 2011.

\bibitem{Raf11b}
James~G. Raftery.
\newblock A perspective on the algebra of logic.
\newblock {\em Quaest. Math.}, 34(3):275--325, 2011.

\bibitem{Ras74}
Helena Rasiowa.
\newblock {\em An algebraic approach to non-classical logics}.
\newblock North-Holland Publishing Co., Amsterdam, 1974.
\newblock Studies in Logic and the Foundations of Mathematics, Vol. 78.

\bibitem{Rau79}
Wolfgang Rautenberg.
\newblock {\em Klassische und nichtklassische {A}ussagenlogik}, volume~22 of
  {\em Logik und Grundlagen der Mathematik [Logic and Foundations of
  Mathematics]}.
\newblock Friedr. Vieweg \& Sohn, Braunschweig, 1979.

\bibitem{Ryb84}
Vladimir~V. Rybakov.
\newblock Admissible rules for logics containing {${\rm S4.3}$}.
\newblock {\em Sibirsk. Mat. Zh.}, 25(5):141--145, 1984.

\bibitem{Ryb97}
Vladimir~V. Rybakov.
\newblock {\em Admissibility of logical inference rules}, volume 136 of {\em
  Studies in Logic and the Foundations of Mathematics}.
\newblock North-Holland Publishing Co., Amsterdam, 1997.

\bibitem{Skv99}
Dmitrij Skvortsov.
\newblock On {P}rucnal's theorem.
\newblock In {\em Logic at work}, volume~24 of {\em Stud. Fuzziness Soft
  Comput.}, pages 222--228. Physica, Heidelberg, 1999.

\bibitem{Slo12}
Katarzyna S{\l}omczy{\'n}ska.
\newblock Algebraic semantics for the {$(\leftrightarrow,\neg\neg)$}-fragment
  of {IPC}.
\newblock {\em MLQ Math. Log. Q.}, 58(1-2):29--37, 2012.

\bibitem{Slo12b}
Katarzyna S{\l}omczy{\'n}ska.
\newblock Unification and projectivity in {F}regean varieties.
\newblock {\em Log. J. IGPL}, 20(1):73--93, 2012.

\bibitem{Woj73}
Ryszard W{\'o}jcicki.
\newblock Matrix approach in methodology of sentential calculi.
\newblock {\em Studia Logica}, 32:7--39, 1973.

\bibitem{Woj98}
Ryszard W{\'o}jcicki.
\newblock {\em Lectures on propositional calculi}.
\newblock Ossolineum Publishing Co., Wroc\l aw, 1984.
\newblock URL=\url{http://www.ifispan.waw.pl/studialogica/wojcicki/}.

\bibitem{ZW12}
Frank Wolter and Michael Zakharyashchev.
\newblock On the {B}lok-{E}sakia theorem.
\newblock In {\em Trends in Logic}, volume in memory of Leo Esakia. To appear.

\bibitem{Wro95}
Andrzej Wro{\'n}ski.
\newblock Transparent unification problem.
\newblock {\em Rep. Math. Logic}, (29):105--107 (1996), 1995.
\newblock First German-Polish Workshop on Logic \& Logical Philosophy
  (Bachotek, 1995).

\bibitem{Wro09}
Andrzej Wro{\'n}ski.
\newblock Overflow rules and a weakening of structural completeness.
\newblock In Janusz Sytnik-Czetwerty{\'n}ski, editor, {\em Rozwa{\.z}ania o
  {F}ilozfii {P}rawdziwej. Jerzemu Perzanowskiemu w Darze}. Wydawnictwo
  Uniwersytetu Jagiello{\'n}skiego, Krak\'ow, 2009.

\end{thebibliography}

\section{Appendix}

Here we add some proof that are not needed to follow the paper. But they say a bit more about considered algebras. Only the fact about presenting $\4^2$ is new.

\subsection{Defining formulas for Heyting algebra $\2^2$ and closure algebra $\4^2$}

\begin{fact}[cf. Lemma \ref{lem:: 2^2 finitely presented}]
Let $\V$ be a nontrivial variety of Heyting algebras. Then $\2^2$ is isomorphic to the $\V$-finitely presented algebra $\P_{x\mjoin \mneg x\= 1}$.
\end{fact}

\begin{proof}
It is known that the variety $\sf V(\2)$, which is term equivalent to the variety of Boolean algebras, is the only minimal variety of Heyting algebras and is defined relative to $\V$ by the identity $(\forall x)[x\mjoin \mneg x\= 1]$. Moreover, $\2^2$ is free of rank one for $\sf V(\2)$. Thus the statement follows.

There also exists a less abstract argument for this. Namely, one may compute that if $\H$ is a Heyting algebra generated by an element $b$ and in which the equality $b\mjoin \mneg b=1$ holds, then the set $\{0,1,b,\mneg b\}$ is closed under basic operations, and hence equals $H$. (The less trivial part of this computation is the verification that $\mneg\mneg b=b$.)
\end{proof}

\begin{fact}[cf. Lemma \ref{lem:: 4^2 is finitely presented}]
Let $\V$ be a variety of closure algebras containing $\4$. Then the algebra $\4^2$ is a $\V$-finitely presented algebra and is isomorphic to  $\P_{(\mBox\Diamond\mBox v\=\Diamond\mBox v) \meet (\Diamond\mBox v \mmeet v\= \mBox v) \meet (\Diamond\mBox v \mjoin v\= \Diamond v)}$.
\end{fact}

\begin{proof}
Let $\alpha$ be the congruence of a free algebra $\Fr_\V(1)$ over $\{v\}$ generated by the pairs
\begin{align*}
e_0=&(\mBox\Diamond\mBox v,\Diamond\mBox v),\\
e_1=&(\Diamond\mBox v \mmeet v, \mBox v),\\
e_2=&(\Diamond\mBox v \mjoin v, \Diamond v).
\end{align*}
The fact that $e_0\in\alpha$ guarantees that $\Diamond\mBox v/\alpha$ is not just closed, but also an open element in $\Fr_\V(1)/\alpha$. Thus $\Fr_\V(1)/\alpha$ is isomorphic to a product $\M_0\times \M_1$ and $\Diamond\mBox a_0=1$ in $\M_0$, $\Diamond\mBox a_1=0$ in $\M_1$, where $a_0$  and $a_1$ are generators of $\M_0$ and $\M_1$ obtained by projecting $v/\alpha$. Now $e_1\in\alpha$ yields that in $\M_0$
\begin{align*}
\mBox a_0&=\Diamond\mBox a_0\mmeet a_0=a_0,\\
\Diamond a_0&=\Diamond\mBox a_0=1
\end{align*}
and hence
\begin{align*}
\mBox \mneg a_0&=0,\\
\Diamond \mneg a_0&=\mneg a_0.
\end{align*}
In particular, $M_0=\{0,1,a_0\mneg a_0\}$ and $\M_0$ is a homomorphic image of $\4$.
No we compute in $\M_1$. Since $e_1\in\alpha$
\[
\mBox a_1=\Diamond\mBox a_1\mmeet a_1=0,
\]
and since $e_2\in\alpha$
\[
\Diamond a_1=\Diamond\mBox a_1\mjoin a_1=a_1.
\]
Hence
\begin{align*}
\mBox \mneg a_1&=a_1,\\
\Diamond \mneg a_1&=0
\end{align*}
and  $\M_1$ is also a homomorphic image of $\4$.  Thus $\Fr_\V(1)/\alpha$ is a homomorphic image of $\4^2$. Now let $h\colon\Fr_\U(1)\to\4^2$ be a surjective homomorphism such that $h(v)=(a,\mneg a)$, where $a$ is the open  and $\mneg a$ is the closed atom in $\4$. A routine verification reveals  that $e_0,e_1,e_2\in\ker(h)$. Thus $\Fr_\V(1)/\alpha$ is isomorphic to $\4^2$.

\end{proof}

\subsection{Splittings for varieties McKinsey algebras and monadic algebras}

\begin{fact}[Proposition \ref{prop:: splitting McKinsey} \protect{\cite[Example III.3.9]{Blo76}\cite[Example IV.5.4]{Rau79}}]
Let $\U$ be a variety of closure algebras. Then
$\S_2\not\in\U$  if and only if  $\U$ satisfies McKinsey identity.
\end{fact}

\begin{proof}
Since $\S_2$ does not satisfy McKinsey identity, the backward implication holds.

For the foreword implication assume that $\U$ does not satisfy McKinsey identity. All modal algebras satisfies the identity
\[
(\forall x)\;[\mu(x)\=\mneg \mBox (\Diamond x   \mmeet  \Diamond \mneg x )].
\]
Hence a modal algebra is a McKinsey algebra iff it satisfies
\[
(\forall x)\;\mBox (\Diamond x   \mmeet  \Diamond \mneg x )\=0.
\]
Therefore there is an algebra $\M\in \U$ and an element $a\in M$ such that
\[
\mBox (\Diamond a   \mmeet  \Diamond \mneg a )>0.
\]
Because $\M$ is a closure algebra, the algebra $\M/\theta$ is nontrivial, where $\theta$ is the congruence of $\M$ generated by $(\Diamond a   \mmeet  \Diamond \mneg a,1)$. Let $\N$ be the subalgebra of $\M/\theta$ generated by $b=a/\theta$. Since  $\Diamond b=\Diamond \mneg b=1$, elements 0 and 1 are the only closed elements in $\N$, and thus $\N$ is isomorphic to $\S_2$ or to $\2$. But this equations also yields that $0<b<1$. Hence $\N$ must be isomorphic to $\S_2$.
\end{proof}

The next splitting is probably the most commonly known, see e.g. \cite[Example on p. 336]{Kra99}  

\begin{fact}
For every variety $\V$ of closure algebras, $\4\not\in \V$ if and only if $\V$ is a variety of monadic algebras. In particular, if $\mathcal Y$ is a non-minimal variety of Heyting algebras, then $\4\in\sigma(\mathcal Y)$.
\end{fact}

\begin{proof}
Recall that the variety of monadic algebras is defined relative to the variety of closure algebras by the identity $(\forall x)[\Diamond\mBox x=\mBox x]$. As $\4$ is not monadic, one direction is clear. For the other direction assume that there is an algebra $\M\in \V$ and an element $a\in M$ such that $\Diamond\mBox a>\mBox a$. Note that a closure algebra $\N$ has a subalgebra isomorphic to $\4$ iff it has an element $d$ such that $d<\Diamond d=1$. If $\mBox\Diamond\mBox a>\mBox a$, then for $\N$ we may take $\M/\theta$, where $\theta$ is the congruence generated by the pair $(\mBox\Diamond\mBox a,1)$, and as $d$ the element $\mBox a/\theta$. So let us assume that $\mBox\Diamond\mBox a=\mBox a$. This case is more difficult, but one may take the algebra ${\sf B}(\2^2\oplus\1)$ from Figure \ref{fig::B(Lev_2^+)} as a prototypical example (and actually only this example is needed in Theorem \ref{thm:: modal companions of Medvedev and Levin}). Indeed, in ${\sf B}(\2^2\oplus\1)$ we have $\Diamond\mBox a>\mBox a=a$  and $a=\mBox a=\mBox\Diamond\mBox a$. However the element $a\mjoin b$ generates a subalgebra isomorphic to $\4$. In general, in $\M$ we may take
\[
d=\mneg\Diamond\mBox a\mjoin \mBox a.
\]
Indeed, Since $\Diamond\mBox a>\mBox a$, we have $d<1$. Moreover,
\[
\Diamond d=\Diamond\mneg\Diamond \mBox a\mjoin \Diamond\mBox a=
\mneg\mBox\Diamond \mBox a\mjoin \Diamond\mBox a=\mneg\mBox a\mjoin \Diamond\mBox a \geq \mneg\mBox a\mjoin \mBox a=1.
\]
\end{proof}

\end{document}